\documentclass[12pt, a4paper]{article}

\usepackage{arxiv}
\usepackage{amsmath}
\usepackage{amsthm}
\usepackage{amssymb}
\usepackage{color}
\usepackage{pdfpages}
\usepackage[utf8]{inputenc} 
\usepackage[T1]{fontenc}    
\usepackage[hidelinks]{hyperref}       
\usepackage{url}            
\usepackage{booktabs}       
\usepackage{amsfonts}       
\usepackage{nicefrac}       
\usepackage{microtype}      
\usepackage{graphicx}
\usepackage{lineno}

\newtheorem{theorem}{Theorem}
\newtheorem{proposition}[theorem]{Proposition}
\newtheorem{lemma}[theorem]{Lemma}

\newtheorem{corollary}[theorem]{Corollary}
\newtheorem{conjecture}[theorem]{Conjecture}

\newtheorem{remark}[theorem]{Remark}

\def\eu{{\ensuremath{\mathrm{e}}}}

\newcommand{\Rset}{\mathbb{R}}
\newcommand{\Cset}{\mathbb{C}}

\title{On some inequalities for the two-parameter Mittag-Leffler function in the complex plane\footnotemark}
\author{Roberto Garrappa \\
Department of Mathematics \\ Università degli Studi di Bari, Italy\\
\texttt{roberto.garrappa@uniba.it}
\And
Stefan Gerhold \\
Institute of Statistics and Mathematical Methods in Economics \\
TU Wien, Austria \\
\texttt{sgerhold@fam.tuwien.ac.at}
\And
Marina Popolizio\\
Department of Electrical and Information Engineering \\
Politecnico di Bari, Italy\\
\texttt{marina.popolizio@poliba.it}
\And
Thomas Simon\\
Laboratoire Paul Painlev{\'e} \\
Universit{\'e} de Lille, France\\
\texttt{thomas.simon@univ-lille.fr}
}

\renewcommand{\headeright}{}
\renewcommand{\undertitle}{}
\renewcommand{\shorttitle}{Some inequalities for the two-parameter Mittag-Leffler function}

\hypersetup{
pdftitle={On some inequalities for the two-parameter Mittag-Leffler function in the complex plane},
pdfauthor={R.~Garrappa, S.~Gerhold, M.~Popolizio, T.~Simon},
pdfkeywords={Mittag-Leffler functions, inequalities, asymptotic expansions},
}

\numberwithin{equation}{section}

\begin{document}

\maketitle


\renewcommand{\thefootnote}{\fnsymbol{footnote}}
\footnotetext[1]{This is the preprint of the paper published in Journal of Mathematical Analysis and Applications,
551(1) (2025), 129588, with doi: \href{https://doi.org/10.1016/j.jmaa.2025.129588}{10.1016/j.jmaa.2025.129588}}

\begin{abstract}
For the two-parameter Mittag-Leffler function $E_{\alpha,\beta}$ with $\alpha > 0$ and $\beta \ge 0,$ we consider the question whether $|E_{\alpha,\beta}(z)|$ and  $E_{\alpha,\beta}(\Re z)$ are comparable on the whole complex plane. We show that the inequality $|E_{\alpha,\beta}(z)|\le E_{\alpha,\beta}(\Re z)$ holds globally if and only if $E_{\alpha,\beta}(-x)$ is completely monotone on $(0,\infty)$. For $\alpha\in [1,2)$ we prove that the complete monotonicity of $1/E_{\alpha,\beta}(x)$ on $(0,\infty)$ is necessary for the global inequality $|E_{\alpha,\beta}(z)|\ge E_{\alpha,\beta}(\Re z),$ and also sufficient for $\alpha =1.$ For $\alpha \ge 2$ we show that the absence of non-real zeros for $E_{\alpha,\beta}$ is sufficient for the global inequality $|E_{\alpha,\beta}(z)|\ge E_{\alpha,\beta}(\Re z),$ and also necessary for $\alpha =2.$ All these results have an explicit description in terms of the values of the parameters $\alpha,\beta.$ Along the way, several inequalities for $E_{\alpha,\beta}$ on the half-plane $\{\Re z \ge 0\}$ are established, and a characterization of its log-convexity and log-concavity on the positive half-line is obtained.
\end{abstract}



\section{Introduction}

In this paper we study several inequalities of the two-parameter Mittag-Leffler function $E_{\alpha,\beta}$, which is defined on the complex plane for $\alpha > 0$ and $\beta\in\Rset$ by the series representation  
\begin{equation}\label{eq:ML}
	E_{\alpha,\beta}(z) = \sum_{k=0}^{\infty} \frac{z^k}{\Gamma(\alpha k + \beta )}
	, \qquad z \in \Cset.
\end{equation}
Mittag-Leffler (ML) functions have a long history deeply intertwined with the development of complex analysis and special functions. One of the most notable recent aspects of ML functions is their pivotal role in fractional calculus. These functions are indeed eigenfunctions of fractional derivatives, and hence they provide analytical tools for studying fractional differential equations -- see \cite{Diethelm2010,KilbasSrivastavaTrujillo2006,Mainardi2022} for more on this topic. They also serve as the basis for devising numerical methods with improved stability properties, see \cite{GarrappaPopolizio2011_CAMWA,SarumiFuratiKhaliqMustapha2021}. 
We refer to Chapter 18 in the classical monograph~\cite{EMOT} and to the more recent
textbook~\cite{GorenfloKilbasMainardiRogosin2020} for an overview of properties and results concerning ML functions.

For $\alpha=\beta=1$ the ML function reduces to the exponential function, namely
\[
E_{1,1}(z) = \eu^z ,
\] 
and for this reason ML functions are sometimes interpreted -- see e.g. Chapters~1 and~2
in~\cite{GorenfloKilbasMainardiRogosin2020} -- as a parametric generalization of the exponential function. Despite this relationship, ML functions do not possess many of the nice properties that characterize the exponential function. For example, the semi-group property $E_{\alpha,\beta}(z+s)\, =\, E_{\alpha,\beta}(z)E_{\alpha,\beta}(s)$ 
is valid only in the exponential case $\alpha=\beta=1$ or for the values $z=s=0$ when $\beta\in\{0,1,2\}$. Instead, one usually obtains some inequalities, and for $\beta = 1,$ it follows from Proposition 2.3 in \cite{GerSim} that
$E_{\alpha,1}(z+s) \ge E_{\alpha,1}(z)E_{\alpha,1}(s)$ for all $\alpha\le 1$ and $z,s\in {\mathbb R}$ while $E_{\alpha,1}(z+s) \le E_{\alpha,1}(z)E_{\alpha,1}(s)$ for all $\alpha\ge 1$ and $z,s\in {\mathbb R}^+.$ A complete characterization of these sub-additivity resp. super-additivity properties for all $\alpha > 0, \beta \ge 0$ will also be given in Remark \ref{Subsuper} (b) below. 

In this paper, we are mainly interested in the following important property of the exponential function:
\begin{equation}\label{eq:Exp_Inequality}
	| \eu^z | = \eu^{\Re z},
\end{equation}
where $\Re z$ denotes the real part of a complex number $z$. The property (\ref{eq:Exp_Inequality}) has, for example, useful applications in studying the behavior of solutions of ordinary differential equations, as well as in analyzing the sensitivity of exponential matrices~$\eu^{A}$ with respect to perturbations of the matrix argument -- see \cite{VanLoan1977}. Investigating some analogue of \eqref{eq:Exp_Inequality} for ML functions is of interest, since ML functions essentially replace the exponential function when moving from integer order to fractional order calculus. In particular, exploring possible extensions of~\eqref{eq:Exp_Inequality} improves knowledge of the sensitivity of matrix ML functions $E_{\alpha,\beta}(A)$ with respect to perturbations of the matrix argument $A$.

As expected, the equality $\bigl| E_{\alpha,\beta}(z)\bigr| = E_{\alpha,\beta}(\Re z)$ cannot hold globally except for $\alpha = \beta = 1,$ as a consequence of our corollary \ref{IFF2}. In this paper, we wish to identify regions in the $(\alpha,\beta)$ parameter space where the inequalities corresponding to \eqref{eq:Exp_Inequality} hold globally, that is, for any $z\in\Cset$. We characterize the global inequality 
\begin{equation}\label{eq:Ineq_LE}
	\bigl| E_{\alpha,\beta}(z)\bigr|\, \le \, E_{\alpha,\beta}(\Re z) 
\end{equation}
by $0<\alpha\le 1$ and $\beta \ge \alpha$, a region which corresponds to the complete monotonicity of the function $E_{\alpha,\beta}(-x)$ on the positive half-line. The characterization of the reverse global inequality  
\begin{equation}\label{eq:Ineq_GE}
	\bigl| E_{\alpha,\beta}(z)\bigr|\, \ge \, E_{\alpha,\beta}(\Re z)
\end{equation}
turns out however to be more subtle. We obtain characterizations on the lines $\{\alpha = 1\}$ with $\beta \le 1,$ and $\{\alpha = 2\}$ with $\beta \le 3,$ and in the remaining regions $\alpha\in (1,2)$, resp. $\alpha > 2$, we obtain a necessary, resp. a sufficient, condition. In the case $\alpha \in (1,2)$ we prove that the global inequality \eqref{eq:Ineq_GE} holds only for $\beta\in [\alpha-1,\alpha]$ and that, in this region, the reciprocal function $1/E_{\alpha,\beta}(x)$ is completely monotonic on $(0,\infty).$ The latter result, which is proved by a specific probabilistic method and has an independent interest, also ensures the validity of \eqref{eq:Ineq_GE} on the right half-plane $\{\Re z \ge 0\}$ for the rhomboid region $\{\alpha\in [1,2], \beta\in[\alpha-1,\alpha]\}.$ We conjecture that \eqref{eq:Ineq_GE} can be extended on the whole complex plane, but contrary to the case $\alpha\le 1$ this extension cannot be obtained directly from the complete monotonicity, since $E_{\alpha,\beta}$ has then negative zeros. 

In the case $\alpha > 2,$ it is known that $E_{\alpha, \beta}$ has an infinite number of negative zeros and we can show by the Hadamard factorization that \eqref{eq:Ineq_GE} holds globally as soon as there are only negative zeros, which happens for $0<\beta \le 2\alpha-1$ or $\alpha \ge 4$ and $\beta \le 2\alpha.$ However, the possibility of the existence of non-real zeros makes the global inequality \eqref{eq:Ineq_GE} less plausible, and we identify a convex unbounded region in the $(\alpha,\beta)$ parameter space where \eqref{eq:Ineq_GE} fails globally. Contrary to the case $\alpha\in (1,2),$ we also conjecture that for $\alpha > 2$ the absence of non-real zeros is actually equivalent to the global inequality \eqref{eq:Ineq_GE}. All in all, it seems to us that a better understanding of the location of the negative zeros of $E_{\alpha,\beta},$ a task which is known to be difficult, however, is necessary to characterize \eqref{eq:Ineq_GE} for $\alpha\not\in\{1,2\}.$ 

\section{Global inequalities}\label{S:MainResults}

In this section we present our main results concerning the extension of \eqref{eq:Exp_Inequality} to Mittag-Leffler functions. We first identify two separate regions for~$\alpha$ and~$\beta$ where $\bigl| E_{\alpha,\beta}(z)\bigr|$ is comparable to $E_{\alpha,\beta}(\Re z)$ on the whole complex plane.

\medskip

\begin{theorem}\label{thm:Ineq_LE}
{\em (a)} Suppose $0<\alpha\le1$ and $\beta \ge \alpha$. Then, the inequality \eqref{eq:Ineq_LE} holds for all $z\in\Cset.$ 

{\em (b)} Suppose $\alpha\ge  2$ and $0\le \beta \le 2\alpha-1$ or $\alpha\ge 4$ and $0\le\beta\le 2\alpha$. Then, the inequality \eqref{eq:Ineq_GE} holds for all $z \in \Cset.$
\end{theorem}

\begin{proof}
We begin with (a). It follows from the main result of \cite{Sc96} that the function $E_{\alpha,\beta}(-x)$ is completely monotone on $(0,\infty)$ if and only if $0<\alpha\leq1$ and $\beta \geq \alpha$ -- see also~\cite{MiSa97} for a very short proof of the ``if'' direction.
Thus, for those $\alpha$ and $\beta$
  there exists a non-negative random variable~$M_{\alpha,\beta}$ such that
  
  \begin{equation}\label{eq:M}
    \Gamma(\beta) E_{\alpha,\beta}(-x) = \mathbb{E}[\eu^{-xM_{\alpha,\beta}}],\qquad x>0.
  \end{equation}
Notice that for $\alpha=1$ and $\beta>1$ the random variable $M_{1,\beta}$ is beta distributed
with parameters $(1,\beta-1)$, whereas for $\alpha\in {(0,1]}$ and $\beta =1$, the random variable $M_{\alpha,1}$ has a so-called Mittag-Leffler distribution -- see \cite{Huillet} and the references therein. Since $E_{\alpha,\beta}$ is an entire function, the analytic continuation theorem
extends~\eqref{eq:M}  to all complex~$z$, whence

$$\Gamma(\beta)|E_{\alpha,\beta}(z)| \,=\, \big|\mathbb{E}[\eu^{zM_{\alpha,\beta}}]\big|
    \,\leq\, \mathbb{E}\big[|\eu^{zM_{\alpha,\beta}}|\big] \,=\, \mathbb{E}\big[\eu^{\Re z M_{\alpha,\beta}}\big]\, =\, \Gamma(\beta)E_{\alpha,\beta}(\Re z)$$
    for all $z\in\Cset.$ 

    \medskip

    We next proceed to (b). It follows from Theorem 3.1.1 in \cite{PopovSedletski2011} and Corollary 2 in \cite{OsPe97} that the zeros of $E_{\alpha,\beta}(z)$ are simple, real and negative for $\alpha> 2$ and $0<\beta \le 2\alpha-1,$ for $\alpha\ge 4$ and $0<\beta\le 2\alpha$ as well as for $\alpha=2$ and $0<\beta<3.$ This property remains true for $\alpha = 2$ and $\beta = 3$ since
    
$$E_{2,3}(z) = \frac{\cosh\sqrt{z}-1}{z}\, = \,\frac{2}{z} \sinh^2 \frac{\sqrt{z}}{2}\, =\,  \frac{1}{2} \prod_{n=1}^{\infty} \Bigl( 1 + \frac{z}{4 n^2 \pi^2} \Bigr)^2$$
by Euler's product formula for sines. Since $E_{\alpha,\beta}(z)$ has finite order $1/\alpha<1$ as an entire function, the Hadamard factorization theorem -- see e.g.\ \cite[Theorem XI.3.4]{Co78} -- implies 

 $$   \Gamma(\beta) E_{\alpha,\beta}(z) = \prod_{n=1}^\infty\Big(1+\frac{z}{a_{n}}\Big),\quad
    z\in\mathbb{C},$$
  where $ 0 < a_{1} < a_{2} < \cdots$ are the absolute values of the zeros of $E_{\alpha,\beta}$. The assertion for $\alpha\ge  2$ and $0 < \beta \le 2\alpha-1$ is a consequence of 
\begin{equation*}
    |1+z/a_n|^2 \,= \,(1+\Re z/a_n)^2+(\Im z/a_n)^2\, \geq \, (1+\Re z/a_n)^2. 
\end{equation*}  
In the remaining case $\alpha \ge 2$ and $\beta = 0,$ we observe that $E_{\alpha, 0} (z) = z E_{\alpha, \alpha} (z)$, so that the statement follows from the case $\alpha = \beta \ge 2$ which was dealt with  previously. 
\end{proof}

\begin{remark}
\label{Rk1}
(a) In the special case $\alpha= 1$ and $\beta \ge 1$ we have that 

$$E_{1,\beta}(z)\, =\, \frac{1}{\Gamma(\beta)} \sum_{n\ge 0} \frac{(1)_n}{(\beta)_n} \frac{z^n}{n!}\, =\, \frac{1}{\Gamma(\beta)}\, {\rm M} (1,\beta, z)$$
is a Kummer confluent hypergeometric function, which in the general case $\beta > \alpha > 0$ admits the integral representation
\[
{\rm M}(\alpha,\beta,z) \, =\, \frac{\Gamma(\beta)}{\Gamma(\alpha)\Gamma(\beta-\alpha) }
\int_0^1 t^{\alpha-1}(1-t)^{\beta-\alpha-1} \eu^{tz}\, dt .
\]

In general, observe that for all $\alpha, \beta > 0$ one has
\[
  \big| {\rm M}(\alpha,\alpha+\beta,z)\big| \, \le\, {\rm M}(\alpha,\alpha+\beta,\Re z), \quad z\in\Cset,
\]
since $M(\alpha,\alpha+\beta,z)$ is then the moment generating function of a beta distribution with parameters $(\alpha,\beta).$ 

\medskip

(b) It remains an open question to determine the exact set of positive parameters $(\alpha,\beta)$ where the zeros of $E_{\alpha,\beta}(z)$ are all negative and simple. We refer to the end of Section 4.6 in \cite{GorenfloKilbasMainardiRogosin2020} and the references therein for more on this topic. 
\end{remark}

Later on, we will see that $\{\alpha \in {(0,1]}, \, \beta\ge \alpha\}$ is the {\em only} set of positive parameters where \eqref{eq:Ineq_LE} holds globally. From the above proof, one might be tempted to believe that the global inequality \eqref{eq:Ineq_GE} is characterized by the negativity and simplicity of all zeros, but we will prove that this is not the case. We begin with a complete monotonicity result which is of an independent interest, and which will be proved by a specific probabilistic argument (it would have been nice to obtain this property by purely analytical techniques). We recall that a smooth function $f : (0,\infty) \to (0,\infty)$ is completely monotone (CM) if $(-1)^n f^{(n)}(x)\ge 0$ for all $n\ge 0$ and $x > 0.$

\medskip

\begin{proposition}
\label{CM}
For every $\alpha\in [1,2]$ and $\beta \in [\alpha-1, \alpha],$ the function 

$$x\,\mapsto\, \frac{1}{E_{\alpha,\beta}(x)}$$ 
is {\em CM} on $(0,\infty).$
\end{proposition}

\begin{proof} The argument relies on certain properties of positive self-similar Markov processes with no positive jumps, their recurrent extensions and their first hitting times for points. We first handle the case $\alpha = \beta \in [1,2],$ which is originally due to Patie -- see Corollary 4.1 in \cite{Patie09}, and also Theorem~2.5 and Example 1.2.4 in \cite{PatieBartho} -- but we will give details for completeness. Observe first that the cases $\alpha = 1,2$ are plain with $E_{1,1} = \eu^z$ and

$$E_{2,2}(z) \, =\, \frac{\sinh \sqrt{z}}{\sqrt{z}}\, =\, \prod_{n\ge 1} \left( 1 + \frac{z}{n^2\pi^2}\right),$$
so that we can focus on the case $\alpha \in (1,2).$ We follow the lines of Theorem 2.1. in \cite{Patie09}. The formula

$$\frac{\Gamma (\alpha z + \alpha)}{\Gamma(\alpha z + 1)}\, =\, \Gamma (\alpha) \, -\, \frac{1}{\alpha\,\Gamma (1-\alpha)}\int_0^\infty (1-\eu^{-z t}) \frac{\eu^{-t}}{(1-\eu^{-t/\alpha})^\alpha}\, dt,$$
which can be obtained by an evaluation of the standard beta integral and an integration by parts, leads after a further integration by parts to

$$\frac{\Gamma (\alpha z + \alpha)}{\Gamma(\alpha z)}\, =\, \Gamma (\alpha+1)z \, +\, \int_{-\infty}^0 (\eu^{z t} - 1 - zt)\, \varphi_\alpha (t) \, dt\, =\Psi_\alpha (z)$$
for all $z\in\Cset$ with $\Re z > -1,$ where 

$$\varphi_\alpha (t) \, =\, \frac{\alpha -1}{\Gamma(2-\alpha)}\frac{\eu^{-\vert t\vert}}{(1-\eu^{-\vert t\vert /\alpha})^{\alpha+1}},\qquad t \in (-\infty, 0).$$
The function $\Psi_\alpha$ is hence the Laplace exponent of a spectrally negative Lévy process $(Z_t)_{t\ge 0}$ with drift coefficient $\Gamma(\alpha +1) >0,$ in other words we have

$${\mathbb E}\big[\eu^{z Z_t}\big]\, =\, \eu^{t \Psi_\alpha (z)}, \qquad \Re z>-1\;\mbox{and}\;  t\ge 0.$$
According to \cite{Lamperti}, we may associate to this Lévy process
$(Z_t)_{t\ge 0}$ a 1-self-similar Markov process $(X_t)_{t\ge 0}$ taking values in $(0,\infty)$, whose infinitesimal generator

$${\mathcal L}_\alpha f (x) \, =\, \Gamma (\alpha+1) f'(x) \, +\, \frac{1}{x}\int_{-\infty}^0 \big(f(x\eu^t) - f(x) - xf'(x) t\big)\, \varphi_\alpha (t) \, dt$$
acts on ${\mathcal C}^1$ functions $f : {\mathbb R}^+ \to {\mathbb R}.$ It is easy to check that ${\mathcal L}_\alpha x^n \, =\Psi_\alpha (n) x^{n-1}$ for all $n\ge 0.$ Fixing $q > 0$ and setting

$$F_q(z) \, =\, 1 \, +\, \sum_{n\ge 1} \frac{q^n z^n}{\Psi_\alpha(1)\cdots \Psi_\alpha(n)}\, =\, \Gamma(\alpha) E_{\alpha, \alpha} (qz),$$
we deduce ${\mathcal L}_\alpha F_q (x) = q F_q(x)$, which shows that $F_q$ is an eigenfunction for ${\mathcal L}_\alpha,$ in other words that the space-time function $G_q(t,x) = \eu^{-qt} F_q(x)$ is such that $(\partial_t + {\mathcal L}_\alpha) {G_q} = 0.$ The classical Dynkin formula for Markov processes -- see Formula (5.8) p.133 in~\cite{Dynk} -- applied to the space-time process $(t,X_t)$ with infinitesimal generator $\partial_t + {\mathcal L}_\alpha$ and the harmonic function $G_q$ at the bounded stopping time $t\wedge T_a$ implies then
\[
  {\mathbb E}_x \big[ \eu^{-q (t\wedge T_a)} F_q(X_{t\wedge T_a})\big] \, = \, F_q (x)
\]
for all $a > x > 0$ and $t > 0,$ where we have written ${\mathbb P}_x$ for the law of
$(t, X_t)_{t\ge 0}$ starting from $(t= 0, x >0) $ and 
\[
  T_a = \inf\{ t > 0 : X_t = a\},
\]
which is finite a.s.\ since $a > x$ and $X$ drifts to $\infty$ with no positive jumps. Letting $t\to\infty$ and using dominated convergence with $X_{t\wedge T_a} \in [0,a]$ for all $t \ge 0$ leads to
\[
 {\mathbb E}_x \big[ \eu^{-q T_a}\big]  =  \frac{F_q (x)}{F_q(a)}
\]
for all $a > x > 0$ and $q > 0.$ Letting $x\to 0$ and taking $a=1,$ we finally obtain
\[
  \frac{1}{E_{\alpha, \alpha} (q)}\, =\, \Gamma(\alpha)\, {\mathbb E}_0 \big[\eu^{-q T_1}\big],
\]
a CM function as required. The proof for the case $\beta \in {[\alpha -1,\alpha)}$ is analogous, except that we will have to deal with recurrent Markov processes. Except in the obvious case $\alpha = 1$ and $\beta = 0$ with $E_{1,0} (z) = z\eu^z,$ we can write again 
\begin{equation}
    \label{Fq}
    F_q(z) \, =\, \Gamma(\alpha) E_{\alpha, \beta} (qz)\, =\, 1 \, +\, \sum_{n\ge 1} \frac{q^n z^n}{\Psi_{\alpha,\beta}(1)\cdots \Psi_{\alpha,\beta}(n)},
\end{equation}
where 

$$\Psi_{\alpha,\beta} (z) \, =\, \frac{\Gamma(\alpha z +\beta)}{\Gamma (\alpha z + \beta - \alpha)}\, =\, \Psi_\alpha (z + \beta\alpha^{-1} -1)$$
is the Laplace exponent of a spectrally negative Lévy process $(Z_t)_{t\ge 0}$. More precisely, on $\{\Re z > -\beta\alpha^{-1}\}$ we have

$$\Psi_{\alpha,\beta} (z)\, =\, \frac{\Gamma (\beta)}{\Gamma(\beta -\alpha)} \, +\, \Psi'_\alpha (\beta\alpha^{-1} -1) z\, +\, \int_{-\infty}^0 (\eu^{z t} - 1 - zt)\, \eu^{(\beta\alpha^{-1}-1) t}\, \varphi_\alpha (t) \, dt$$
for $\beta \in {(\alpha-1,\alpha)}$ and 

$$\Psi_{\alpha,\alpha-1} (z)\, =\, - \alpha\Gamma (\alpha-1)z \, +\, \int_{-\infty}^0 (\eu^{z t} - 1 - zt)\, \eu^{- \alpha^{-1}t}\, \varphi_\alpha (t) \, dt.$$
The difference with the case $\beta = \alpha$ is that the Lévy process $(Z_t)_{t\ge 0}$ is killed to $-\infty$ at an exponential time with parameter $-\Gamma (\beta)/\Gamma(\beta -\alpha)$ for $\beta \in {(\alpha-1,\alpha)}$ and has a negative drift coefficient $- \alpha\Gamma (\alpha-1)$ for $\beta = \alpha -1.$ As explained in the introduction of \cite{Rivero05}, this shows that the corresponding Lamperti process will have a trap at zero, which is reached in finite time a.s.\ starting from any $x > 0.$ The main result of \cite{Rivero05} makes it however possible to build a recurrent extension of this Lamperti process, which we denote again by $(X_t)_{t\ge 0}$. More precisely, the unique positive root of $\Psi_{\alpha,\beta}$ is $\theta = 1- \beta/\alpha < 1$ so that we can apply Theorem 2 (i) in \cite{Rivero05} -- see also Example 4 therein. For smooth functions $f : {\mathbb R}^+ \to {\mathbb R}$ vanishing at zero, it follows from Proposition 1.1. in \cite{Patie09} and the fact that $\theta < 1$ that the action of the infinitesimal generator of $(X_t)_{t\ge 0}$ is given by 
\[
  {\mathcal L}_{\alpha,\beta} f (x) \, =\, \frac{\Gamma (\beta)}{\Gamma(\beta -\alpha)} \, \frac{f(x)}{x}\, +\, \Psi'_\alpha (\beta\alpha^{-1} -1)f'(x) \, +\, \frac{1}{x}\int_{-\infty}^0 \big(f(x\eu^t) - f(x) - xf'(x) t\big)\, \eu^{(\beta\alpha^{-1}-1) t}\,\varphi_\alpha (t) \, dt
\]
for $\beta \in {(\alpha-1,\alpha)}$ and 
\[
  {\mathcal L}_{\alpha,\alpha-1} f (x) \, =\,- \alpha\Gamma (\alpha-1) f'(x) \, +\, \frac{1}{x}\int_{-\infty}^0 \big(f(x\eu^t) - f(x) - xf'(x) t\big)\, \eu^{- \alpha^{-1}t}\,  \varphi_\alpha (t) \, dt
\]
for $\beta = \alpha -1$. We hence again obtain ${\mathcal L}_{\alpha, \beta} x^n \, =\Psi_{\alpha,\beta} (n) x^{n-1}$ for all $n\ge 1.$ The action of the infinitesimal generator for smooth functions non-vanishing at zero is more complicated and depends on the construction of the recurrent extension -- 
see Section~2 in~\cite{Rivero05}. However, the conservativeness of the recurrent
extension  $(X_t)_{t\ge 0},$ which means $P_t 1 = 1$ where $(P_t)_{t\ge 0}$ is the associated semi-group, implies readily ${\mathcal L}_{\alpha,\beta} 1 = 0$ by definition of the infinitesimal generator, and we again obtain
\[
  {\mathcal L}_{\alpha,\beta} F_q\, =\, q F_q,
\]
where $F_q$ is the function defined in~\eqref{Fq}. Reasoning exactly as above yields
\[
  \frac{1}{E_{\alpha, \beta} (q)}\, =\, \Gamma(\beta)\, {\mathbb E}_0 \big[\eu^{-q T_1}\big], \qquad q\ge 0,
\]
a CM function as required, where $T_1 = \inf\{t > 0 : X_t =1\}$ is a.s. finite under~${\mathbb P}_0$, since~$X$ is recurrent, which implies that $S_1 = \inf\{t > 0 : X_t >1\}$ is a.s. finite, and does not have positive jumps, which yields $T_1 = S_1$ a.s. 
\end{proof}

\begin{remark} 
\label{CBFFF}
(a) It follows easily from the absence of positive jumps of the Markov process $(X_t)_{\ge 0}$ that the underlying random variable $T_1$ is $n$-divisible for every $n\ge 1,$ and is hence infinitely divisible. This shows that there is a Bernstein function $\varphi_{\alpha,\beta}(z)$ such that $E_{\alpha, \beta}$ admits the exponential representation 

$$E_{\alpha, \beta} (z)\, =\, \frac{\eu^{\varphi_{\alpha,\beta}(z)}}{\Gamma(\beta)}$$
on $\{\Re z \ge 0\}$ for all $\alpha\in [1,2]$ and $0\neq \beta \in [\alpha-1,\alpha].$ Notice that except for $\alpha = 2$ and $(\alpha,\beta) =(1,1)$ the function $E_{\alpha,\beta}$ has an infinite number of non-real zeros so that $\varphi_{\alpha,\beta}$ is not a complete Bernstein function. The non-complete character of the function $\varphi_{\alpha,\beta}$ will be discussed in further detail in Section \ref{S:Conj}. We refer to \cite{SSV} for a classic account on Bernstein functions. 

\medskip

(b) By the aforementioned Hadamard factorization, the function $1/E_{\alpha,\beta}(x)$ is {\em CM} on $(0,\infty)$ as soon as all the zeros of $E_{\alpha,\beta}$ are negative, in particular also for all $\alpha\ge 2$ and $\beta \in [0,2\alpha-1]$ and for all $\alpha\ge 4$ and $\beta \in [0,2\alpha].$ In the case $\beta =1,$ the complete monotonicity of the function $1/E_{\alpha} (x)$ for all $\alpha\ge 1$ had been stated as an open problem in \cite{GerSim}, where it was shown by other methods that the function is log-convex. See Proposition 2.3 therein and the comments thereafter. It follows from Theorem \ref{thm:Ineq_LE} and from Proposition \ref{prop:Ineq_large_LE} below that for all $\alpha \le 1$ and $\beta \ge 0,$ one has the characterization 

$$\mbox{{\em The function $1/E_{\alpha,\beta}$ is CM on $(0,\infty)$}}\quad\Longleftrightarrow\quad \alpha =1\;\mbox{{\em and}}\;\beta \le 1.$$ 

It is an interesting open problem to characterize the {\em CM} property of $1/E_{\alpha,\beta}$ on $(0,\infty)$ for all $\alpha > 1$ and $\beta \ge 0.$  See below Remark~\ref{b34}~{\em (b)} for some partial observations, and Conjecture  \ref{conj:rhombus} which states that in the case $\alpha\in (1,2],$ the above criterion $\beta\in[\alpha-1, \alpha]$ is indeed a characterization. 
\end{remark}

The following consequence of Proposition~\ref{CM} shows the aforementioned fact that there exist ML functions $E_{\alpha,\beta}$ having an infinite number of complex zeros and nevertheless satisfying \eqref{eq:Ineq_GE} on the whole complex plane.

\medskip

\begin{corollary}\label{CHF}
For all  $\beta\in [0,1],$ the function $E_{1,\beta}(z)$ satisfies \eqref{eq:Ineq_GE} for all $z\in\Cset.$
\end{corollary}

\begin{proof}
We discard the trivial situations $\beta = 0$ with $E_{1,0} (z) = z\eu^z$ and $\beta = 1$ with $E_{1,1}(z) = \eu^z.$ In the remaining case $\beta\in {(0,1)}$, Proposition~\ref{CM} combined with Bernstein's theorem shows that 
\begin{equation}
\label{CMCHF}
 \frac{1}{E_{1, \beta} (z)}\, =\, \Gamma(\beta)\, {\mathbb E} \big[\eu^{-z T}\big], \qquad \Re z\ge 0, 
\end{equation}
where $T$ is some positive random variable. Moreover, since $\Gamma(\beta)E_{1,\beta}(z) = {\rm M} (1,\beta,z)$ is a Kummer confluent hypergeometric function with $\beta\in {(0,1)}$, it is classical -- see e.g.\ Chapter~6.16 in \cite{EMOT} -- that it has a single real root $x_\beta < 0$. The principle of analytic continuation for Laplace transforms of positive random variables --
see e.g.~\cite{Widder41} -- implies then that~\eqref{CMCHF} holds on $\{\Re z > x_\beta\}.$ Reasoning as in Theorem \ref{thm:Ineq_LE}, we obtain
\[
  \left| \frac{1}{E_{1, \beta} (z)} \right|\, \le\, \frac{1}{E_{1, \beta} (\Re z)}, \qquad \Re z > x_\beta,
\]
which shows that \eqref{eq:Ineq_GE} holds on $\{\Re z > x_\beta\}.$ On the other hand, the asymptotic

$$x {\rm M} (1,\beta,-x)\, \to \, \frac{1}{\Gamma(\beta -1)}\, < \, 0\qquad \mbox{as $x\to\infty,$}$$
which is given e.g. in Formula 6.13.1.(3) in \cite{EMOT}, shows that $E_{1,\beta} (\Re z) \le 0$ on $\{ \Re z \le x_\beta\},$ so that \eqref{eq:Ineq_GE} holds on $\Cset$.
\end{proof}

\begin{remark} In view of the preceding result and Remark~\ref{Rk1}, one may ask if for all $\alpha, \beta > 0$ one has
\begin{equation}
\label{BigM}
\big| M(\alpha,\alpha-\beta,z)\big| \, \ge\, M(\alpha,\alpha-\beta,\Re z), \quad z\in\Cset.
\end{equation}
In this respect, it can be shown by the probabilistic method of Proposition~\ref{CM} that $1/M(\alpha,\alpha-\beta,z)$ is {\em CM} for all 
$\alpha > 1$ and $\beta > 0$, with an underlying spectrally negative Lévy process having Laplace exponent
\[
  \Psi(z)\, =\, \frac{z(z+\alpha-\beta-1)}{z+\alpha-1} \, =\,  \left(\frac{\alpha-1-\beta}{\alpha-1}\right) z\, +\, \beta(\alpha -1)\int_{-\infty}^0 (\eu^{zt} - 1 - zt)\eu^{-(\alpha-1)\vert t\vert} dt.
\]
In particular, \eqref{BigM} holds globally for $\alpha > 1 \ge\beta > 0$ since $M(\alpha, \beta,z)$ has then only one real zero, which is negative -- see again Chapter~6.16 in~\cite{EMOT}. The question whether~\eqref{BigM} holds globally for all $\alpha, \beta > 0$ seems to be still open.  
\end{remark}

The following second consequence of Proposition~\ref{CM} shows that for all $\alpha\in[1,2]$ and $\beta \in {(\alpha-1,\alpha]}$, the inequality~\eqref{eq:Ineq_GE} holds globally, except possibly on some compact subset of the half-plane $\{\Re z < 0\}$.

\medskip

\begin{corollary}\label{thm:Ineq_GE12}
Let $\alpha\in {(1,2)}$ and $\beta \in {(\alpha-1, \alpha]}$. There exists $a\le b < 0 < c$ such that \eqref{eq:Ineq_GE} holds for all $z \in \Cset$ such that $\Re z \notin {(a,b)}$ or $\vert\Im z\vert \ge c$.
\end{corollary}

\begin{proof}
The existence of the lower parameter $a$ is a direct consequence of the known asymptotics of $E_{\alpha,\beta}(z)$ along the negative real axis -- see e.g.\ Formula 18.1.(21) in \cite{EMOT}.  More precisely, if $\alpha\in [1,2)$ and $\beta \in {(\alpha-1,\alpha)}$, we have   

$$z E_{\alpha,\beta}(-z)\, \to \, \frac{1}{\Gamma(\beta -\alpha)}\, < \, 0\qquad \mbox{as $z\to\infty$}$$ so that there exists $a > - \infty$ such that 

$$E_{\alpha,\beta}(\Re z) \, <\, 0 \,\le \, \bigl| E_{\alpha,\beta}(z)\bigr|$$
for all $z \in \Cset$ such that $\Re z \le a.$ If $\beta = \alpha\in {(1,2)}$, the existence of $a$ follows similarly from    

$$z^2 E_{\alpha,\alpha}(-z)\, \to \, -\frac{1}{\Gamma(-\alpha)}\, < \, 0\qquad \mbox{as $z\to\infty.$}$$ 
The parameter $b$ can be defined as the first negative zero of $E_{\alpha,\beta}$ as in Corollary~\ref{CHF}, since by Proposition~\ref{CM} and the principle of analytic continuation, the function

$$z\,\mapsto\, \frac{\Gamma(\beta)}{E_{\alpha,\beta}(z)}$$
is a moment generating function defined on $\{\Re z > b\}$ for all $\alpha\in {(1,2]}$ and $\beta\in[\alpha-1,\alpha],$ so that

$$\left| E_{\alpha, \beta} (z) \right|\, \ge\, E_{\alpha, \beta} (\Re z)$$
on $\{\Re z \ge b\}.$ Finally, the continuous function $E_{\alpha,\beta}(\Re z)$ is bounded on
$\{\Re z\in [a,b]\}$, whereas $\vert E_{\alpha,\beta}(z)\vert \to\infty$ as $\vert \Im z\vert \to\infty$ with $\Re z\in [a,b]$, since then $\vert\arg z\vert \to \pi/2$ and we can apply the estimate (4.4.16) in \cite{GorenfloKilbasMainardiRogosin2020}. This implies the existence of the upper
parameter~$c$. 
\end{proof}

\begin{remark} If $E_{\alpha,\beta}(z)$ has exactly one real zero, then this zero must be negative, and the same argument as in Corollary~\ref{CHF} shows that \eqref{eq:Ineq_GE} holds globally. According to the numerical simulations in~\cite{Hanneken2007}, the function $E_{\alpha,1}(z)$ has only one real zero for $\alpha \in {(1, \alpha_0]}$ with $\alpha_0 = 1.42219069..$ See also \cite{HilfeITSF} for related simulations on the complex conjugate zeros. In general, there does not seem to exist a rigorous analysis of the exact number of real zeros for Mittag-Leffler functions, as is the case for confluent hypergeometric functions. The main hindrance to this rigorous analysis is that, contrary to the Kummer functions, it does not seem possible to build a Sturmian chain associated to Mittag-Leffler functions, which are not solutions to linear {\em ODEs}, but to fractional {\em ODEs}. 
\end{remark}

We would like to conclude this section with the following curious two-sided estimate on the half-plane $\{\Re z \ge 0\}.$

\medskip

\begin{proposition}
\label{2sided}
Suppose $\alpha\in [1,2)$ and $\beta\in[1,\alpha]$. Then, one has
\begin{equation}
\label{eq:2S}
E_{\alpha,\beta}(\Re z)\, \le\,\bigl| E_{\alpha,\beta}(z)\bigr|\, \le\, E_{\alpha,\beta}\big((\Re z^{1/\alpha})^\alpha\big)
\end{equation}
on $\{\Re z \ge 0\}.$ Moreover, if $\alpha = 2$ and $\beta\in [1,3],$ the inequalities \eqref{eq:2S} hold true for all $z\in\Cset\setminus(-\infty,0].$ 
\end{proposition}

\begin{proof} We begin with the case $\alpha \in[1,2).$ The first inequality in \eqref{eq:2S} is a direct consequence of Proposition~\ref{CM} and is actually valid for all $\beta \in [\alpha-1,\alpha].$ In order to show the second inequality, we use the representation
\begin{equation}
\label{Kumm}
E_{\alpha,\beta}(x^\alpha)\, = \, \frac{1}{\alpha\Gamma(\beta)} \left(M(1,\beta,x)\, +\, (\alpha - 1) {\mathbb E}[\eu^{-x X_{\alpha, \beta}}]\right)
\end{equation}
for all $\alpha\in[1,2], \beta \ge 1$ and $x\ge 0$, where $X_{\alpha, \beta}$ is some positive random variable, a representation which is a consequence of Theorem B (c) and (d) in \cite{SimITSF} and Formulas 6.3.(7) and 6.9.2.(22) in \cite{EMOT}. By analytic continuation, we deduce that for all $\alpha\in[1,2)$ and $\beta \ge 1$, on $\{\Re z \ge 0\}\subset\{\Re z^{1/\alpha} \ge 0\}$, one has
\begin{align*}
\bigl| E_{\alpha,\beta}(z)\bigr| & \le  \frac{1}{\alpha \Gamma(\beta)}\,\left(\big|M(1,\beta,z^{1/\alpha})\big|\, +\, (\alpha- 1) \big|{\mathbb E}[\eu^{-z^{1/\alpha} X_{\alpha, \beta}}]\big|\right)\\
& \le \frac{1}{\alpha\Gamma(\beta)}\,\left(M(1,\beta,\Re z^{1/\alpha})\, +\, (\alpha- 1) {\mathbb E}[\eu^{-\Re z^{1/\alpha} X_{\alpha, \beta}}]\right)\, =\, E_{\alpha,\beta}((\Re z^{1/\alpha})^\alpha),
\end{align*}
where in the second inequality we have used Remark \ref{Rk1} (a) and reasoned as in the first part of Theorem~\ref{thm:Ineq_LE}. 

\medskip

We next proceed to the case $\alpha =2,$ where the first inequality in \eqref{eq:2S} is a consequence of Theorem \ref{thm:Ineq_LE} (b). For the second inequality, we use again the representation \eqref{Kumm}, which here simply reads
\[
  E_{2,\beta}(x^2)\, = \, \frac{1}{2 \Gamma(\beta)} \left(M(1,\beta,x)\, +\, M(1,\beta,-x)\right),
\]
as a consequence of Formula 18.1.(24) in \cite{EMOT}, and remains valid on the whole complex plane. Since $\Re(\sqrt{z}) \ge 0$ for every $z\in\Cset\setminus(-\infty,0],$ we deduce
\begin{align*}
\bigl| E_{2,\beta}(z)\bigr| & \le  \frac{1}{2 \Gamma(\beta)}\,\left(\left| M(1,\beta,\sqrt{z})\right|\, +\, \left|M(1,\beta,-\sqrt{z}) \right|\right)\\
& \le  \frac{1}{2\Gamma(\beta)}\,\left(M(1,\beta,\Re\sqrt{z})\, +\, M(1,\beta,-\Re\sqrt{z})\right)\, =\, E_{\alpha,\beta}((\Re \sqrt{z})^2),
\end{align*}
where in the second inequality we have used again Remark \ref{Rk1} (a). 
\end{proof}

\begin{remark} The inequality 

$$E_{\alpha,\beta}(\Re z)\,\le\, E_{\alpha,\beta}((\Re z^{1/\alpha})^\alpha)$$ 
holds on $\{\Re z \ge 0\}$ for all $\alpha \ge 1$ and $\beta\ge 0.$ This follows from the positivity of the coefficients of the entire series defining $E_{\alpha,\beta}$ and the inequality $\Re z \le (\Re z^{1/\alpha})^\alpha,$ which amounts to $\cos(\theta) \, \le \, (\cos(\theta/\alpha))^\alpha$ for $\vert \theta\vert < \pi/2$ and $\alpha\ge 1.$ The latter inequality is an immediate consequence of the log-concavity of the cosine on $({-\pi/2},\pi/2)$. 
\end{remark}

\section{Local inequalities}
\label{S:LocalInequalities}

In this section we analyze the inequalities \eqref{eq:Ineq_LE} and \eqref{eq:Ineq_GE} locally, with the help of the following expansion for fixed $\Re z$ as $\Im z \to 0$:
\begin{equation}
\label{LCV}
| E_{\alpha,\beta}(z)|^2 \, -\, \left(E_{\alpha,\beta}(\Re z)\right)^2 \, =\, F_{\alpha,\beta}(\Re z) \times \left(\Im z\right)^2 \, +\, {\rm O} \left(\Im z\right)^4
\end{equation}
with the notation $F_{\alpha,\beta} = (E_{\alpha,\beta}')^2 - E_{\alpha,\beta}E_{\alpha,\beta}'',$ which follows readily from the asymptotic expansion
\begin{equation}
\label{Expa}
E_{\alpha,\beta} (z) \, \sim\, \sum_{k\ge 0} E_{\alpha, \beta}^{(k)} (\Re z)\, \frac{({\rm i}\, \Im z)^k}{k!}
\end{equation}
as $\Im z\to 0.$ We will use the following estimate.

\medskip

\begin{proposition}
\label{LCX}
With the above notation, one has 

$$F_{\alpha,\beta}(x)\, =\, \frac{x^{\frac{2(1-\alpha-\beta)}{\alpha}}\, \eu^{2x^{1/\alpha}}}{\alpha^4} \left( (\alpha-1)x^{1/\alpha}\, + \, \alpha(1-\beta)\, +\,{\rm o}(1) \right) \qquad\mbox{as $x\to\infty.$}$$ 
\end{proposition}

\begin{proof} Equating the entire series, we have the decompositions

$$E_{\alpha,\beta}'\, =\, \frac{1}{\alpha}\left( E_{\alpha, \alpha+\beta-1} + (1-\beta) E_{\alpha, \alpha+\beta}\right)$$
and
$$E_{\alpha,\beta}''\, =\, \frac{1}{\alpha^2}\left( E_{\alpha, 2\alpha+\beta-2}\, + \, (3-\alpha-2\beta) E_{\alpha, 2\alpha+\beta-1}\, + \, (1-\beta)(1-\alpha-\beta) E_{\alpha, 2\alpha+\beta}\right),$$
which lead to
\begin{align*}
F_{\alpha,\beta}& = \frac{1}{\alpha^2}\Big( \big( E_{\alpha, \alpha+\beta-1} + (1-\beta) E_{\alpha, \alpha+\beta}\big)^2\\
&  \quad \, -\,  E_{\alpha,\beta}\big( E_{\alpha, 2\alpha+\beta-2}\, + \, (3-\alpha-2\beta) E_{\alpha, 2\alpha+\beta-1}\, + \, (1-\beta)(1-\alpha-\beta) E_{\alpha, 2\alpha+\beta}\big)\Big).
\end{align*}
The behaviour of the above function for large values on the positive half-line is then obtained after plugging in the estimate 18.1.(22) in \cite{EMOT} for the corresponding ML functions, whose leading term is given by~$t_0$ with the notation therein, and making some algebraic simplifications. We omit details.
\end{proof}

We can now show the aforementioned fact that the global inequality \eqref{eq:Ineq_LE} only holds for the parameter values appearing in Theorem \ref{thm:Ineq_LE} (a).

\medskip

\begin{corollary}
\label{IFF1}
The inequality \eqref{eq:Ineq_LE} holds for all $z\in\Cset$ if and only if $0<\alpha\le 1$ and $\beta \ge \alpha$.
\end{corollary}

\begin{proof}
The if part was proved in Theorem \ref{thm:Ineq_LE} (a). For the only if part, we first observe from \eqref{LCV} and Proposition \ref{LCX} that if $\alpha > 1$ or $\alpha = 1$ and $\beta < 1,$ then we  have

$$| E_{\alpha,\beta}(z)|\, > \, | E_{\alpha,\beta}(\Re z)|\, = \, E_{\alpha,\beta}(\Re z)$$
for $\Re z$ large enough and $\Im z\neq 0$ small enough. Finally, if $0\le\beta < \alpha <1,$ the estimate 18.1.(21) in \cite{EMOT} implies 

$$-(\Re z)\, E_{\alpha,\beta}(\Re z)\, \to\, \frac{1}{\Gamma(\beta -\alpha)}\, <\, 0\qquad \mbox{as $\,\Re z\,\to\, -\infty,$}$$
so that \eqref{eq:Ineq_LE} cannot hold either.
\end{proof}

Another direct consequence of Proposition \ref{LCX} is the following characterization of the exponential function among two-parameter Mittag-Leffler functions.

\medskip

\begin{corollary}
\label{IFF2}
  Let $\alpha>0$ and $\beta\geq0$. Suppose that the equality $\bigl| E_{\alpha,\beta}(z)\bigr| = E_{\alpha,\beta}(\Re z)$ holds in the angular sector $\vert\arg z\vert <\varepsilon$ for some $\varepsilon > 0$. Then, one has $\alpha=\beta=1.$
\end{corollary}

Proposition \ref{LCX} also implies that the global inequality \eqref{eq:Ineq_GE} only holds for $\alpha > 1$ or for $\alpha = 1$ and $\beta \le 1.$ The case $\alpha =1$ being already settled in Corollary \ref{CHF}, we will focus in the remainder of this section on the case $\alpha > 1.$ We begin with the case $\alpha\in {(1,2)}$ and show that the rhomboid region $\beta\in [\alpha-1,\alpha],$ corresponding to the complete monotonicity result of Proposition~\ref{CM}, is the only one where \eqref{eq:Ineq_GE} can hold globally.

\medskip

\begin{proposition}\label{prop:Ineq_large_LE}
Suppose $\alpha\in {(1,2)}$ and $\beta\in {[0,\alpha-1)}\cup{(\alpha,\infty)}$. Then, the inequality \eqref{eq:Ineq_GE} cannot hold globally.
\end{proposition}

\begin{proof}

Applying again Formula 18.1.(21) in \cite{EMOT} gives 

$$-(\Re z)\, E_{\alpha,\beta}(\Re z)\, \to\, \frac{1}{\Gamma(\beta -\alpha)}\, >\, 0\qquad \mbox{as $\,\Re z\,\to\, -\infty,$}$$ where the inequality comes from the fact that $\alpha\in {(1,2)}$ and $\beta\in [0,\alpha-1)\,\cup\,(\alpha,\infty)$. On the other hand, the same formula 18.1.(21) in \cite{EMOT} combined with the decompositions of $E_{\alpha,\beta}'$ and $E_{\alpha,\beta}''$ which were used during the proof of Proposition \ref{LCX} yield after some algebraic simplifications the estimate

$$-(\Re z)^4\, F_{\alpha,\beta}(\Re z)\, \to\, \frac{1}{\Gamma(\beta -\alpha)^2}\, >\, 0\qquad \mbox{as $\,\Re z\,\to\, -\infty.$}$$
Putting those two estimates together with \eqref{LCV} implies the inequality 

$$| E_{\alpha,\beta}(z)|\, <\, E_{\alpha,\beta}(\Re z)$$ for $\Re z$ close enough to $-\infty$ and $\Im z$ close enough to zero, which completes the argument.
\end{proof}

In order to handle the remaining situation $\alpha \ge 2,$ we will need the following preparatory lemma involving a certain function already implicitly appearing in \cite{PopovSedletski2011} for the study of the existence of non-real zeros for $E_{\alpha,\beta}$ -- see in particular the inequality (3.11.1) therein. This function will also play a role in characterizing the log-convexity and the log-concavity of $E_{\alpha,\beta}$ on the positive half-line -- see Theorem \ref{LCLX} below. 

\medskip

\begin{lemma}\label{lem:h_Function}
For each $x > 0,$ there exists a unique solution $y = h(x) > 0$ to the equation
  \begin{equation}\label{eq:h_Function}
    2\Gamma(x+y)^2=\Gamma(y)\Gamma(2x+y).
  \end{equation}
The function $h: (0,\infty) \to (0,\infty)$ is real-analytic and increasing with $\lim_{x\to 0} h(x) = 0$ and $\lim_{x\to 0} h'(x) = \sqrt{2} -1.$ Moreover, it is strictly convex on $[1,\infty)$ with $h(1) = 1, h'(1)=2$ and $h(x) \sim x^2/\log(2)$ as $x\to\infty.$ 
\end{lemma}
\begin{proof}
Setting

$$F(x,y) \, =\, \log 2 \, + \, 2\log \Gamma(x+y)\, -\, \log \Gamma(y)\, -\, \log \Gamma(2x+y)$$
for all $x,y > 0,$ we have $\partial_y F(x,y)= 2\psi(x+y)-\psi(y)-\psi(2x+y) > 0$ by the strict concavity of the digamma function~$\psi.$ Since for each $x > 0$ one has $F(x,y)\to -\infty$ as $y\to 0$ and, by Stirling's formula, $F(x,y)\to  \log 2 > 0$ as $y\to\infty,$ this shows that the function $h(x)$ is properly defined. Moreover, the function $h$ is real-analytic on $(0,\infty)$  by analyticity on $(0,\infty)$ of the gamma function and the analytic implicit function theorem. Setting $G(x) = F(x, h(x)) = 1$ and equating

$$0 \, =\, G'(x)\, =\, \partial_x F(x,h(x))\, +\, h'(x) \partial_y F(x,h(x)),$$
we see that $h'(x) > 0$ for all $x > 0$ because $\partial_x F(x,y) = 2(\psi(x+y)-\psi(2x+y))<0,$
by the increasing character of the digamma function. This shows that $h(x)$ has a non-negative limit as $x\to 0$ and it is clear by \eqref{eq:h_Function} that this limit must be zero. The asymptotic
\[
  2\,\Gamma\big(\sqrt{2}x\big)^2\,\sim\, \frac{1}{x^2}\,\sim\,  \Gamma\big((\sqrt{2}-1)x\big)\, \Gamma\big((\sqrt{2}+1)x\big)\qquad\mbox{as $x\to 0$}
\]
implies $x^{-1} h(x)\to \sqrt{2} -1$ as $x\to 0,$ so that $h$ extends to a ${\mathcal C}^1$ function on $[0,\infty)$ with $h(0) = 0$ and $h'(0) = \sqrt{2}-1.$ We next show the convexity property. Equating again

$$0 \, =\, G''(x)\, =\, H (x,h(x), h'(x)) \, +\, h''(x)\partial_y F(x,h(x))$$
with $H(x,y,z) = \partial_{xx} F(x,y)\, +\, 2z \partial_{xy} F(x,y)\, +\, z^2 \partial_{yy} F(x,y)$ we compute, making some algebraic simplifications and differentiating Formula 1.7.2.(22) in \cite{EMOT},
\begin{align*}
    H(x,y,z) & =  2(z+1)^2 \psi'(x+y) \, -\, z^2\psi'(y) \, -\, (z+2)^2 \psi'(2x+y)\\
    & =  \int_0^\infty \frac{t \eu^{-ty}}{1-\eu^{-t}} \left(2(z+1)^2 \eu^{-tx} \, -\, z^2\, -\, (z+2)^2 \eu^{-2tx}\right) \, dt\\
    & <  -\int_0^\infty \frac{t \eu^{-ty}}{1-\eu^{-t}} \left( (z+2)^2 \eu^{-2tx}\, -\, (z^2+4z + 2) \eu^{-tx}\right) \, dt\\
    & =  -\frac{(z+2)^2}{2}\int_0^\infty \frac{t \eu^{-t(x+y/2)}}{1-\eu^{-t/2}} \left(1 \, -\, 2\left(1 -\frac{2}{(z+2)^2}\right) \frac{\eu^{-ty/2}}{1+\eu^{-t/2}}\right) \, dt\\
    & <  -\frac{(z+2)^2}{2}\int_0^\infty \frac{t \eu^{-t(x+y/2)}}{1-\eu^{-t/2}} \left(1 \, -\, \frac{2\,\eu^{-ty/2}}{1+\eu^{-t/2}}\right) \, dt.
\end{align*}
Therefore, since $2\eu^{-ty/2} < 1+\eu^{-t/2}$ for all $y\ge 1/2$ and $t > 0,$ we have $H(x,y,z) < 0$ for all $y\ge 1/2.$ Putting everything together, we obtain $h''(x) > 0$ for all $x\in[h^{-1}(1/2), \infty)$, and since it is clear from \eqref{eq:h_Function} that $h(1) =1 > h^{-1}(1/2),$ this shows that $h$ is strictly convex on $[1,\infty)$. In particular, the function $x^{-1}h(x)$ increases to some limit as $x\to\infty$, and if this limit~$\ell$ were finite, then \eqref{eq:h_Function} and Stirling's formula would imply $2(\ell +1) = 2\ell +1,$ a contradiction. Hence, one has $x^{-1}h(x)\to\infty$ as $x\to\infty$, and a further application of Stirling's formula
to~\eqref{eq:h_Function} yields

$$\log 2 \, +\, (2x + 2h(x) -1) \log(1+x/h(x))\, = \, (2x+h(x) -1) \log (1+ 2x/(h(x))\, +\, {\rm o}(1),$$
which easily leads to the required asymptotic $h(x) \sim x^2/\log(2)$ as $x\to\infty$. The computation

$$h'(1)\, =\, \frac{2(\psi(3)-\psi(2))}{2\psi(2) - \psi(1) - \psi(3)}\, =\, 2$$ 
finishes the proof.
\end{proof}

\begin{remark} 
\label{ZZZ}
{\em (a)} We believe that $h$ is strictly convex on the whole $(0,\infty).$ On the other hand, it can be shown from the convexity of the function $1/\psi'$ on $(0,\infty)$ -- see \cite{Alz} -- that the function $H(x,y,z)$ takes positive values on $(0,\infty)^3.$ It seems that a more technical analysis is necessary in order to prove the convexity of $h$ on $(0,h^{-1}(1/2)].$

\medskip

{\em (b)} A further application of Stirling's formula, whose details are skipped, yields

$$h(x)\, = \, x^2/\log(2) - x + (11\log(2)/3 -1) + {\rm o}(1)\qquad\mbox{as $x\to\infty.$}$$
\end{remark}

The following proposition depicts a certain region of the parameter space where \eqref{eq:Ineq_GE} does not hold globally, which is relevant for $\alpha \ge 2.$ 

\medskip

\begin{proposition}\label{prop:Ineq_small}
Let $h$ be the function defined in Lemma~\ref{lem:h_Function}. Then, \eqref{eq:Ineq_GE} does not hold globally if $\beta>h(\alpha).$
\end{proposition}

\begin{proof} With the above notation, one has

$$F_{\alpha, \beta}(0)\, =\, \frac{1}{\Gamma(\alpha+\beta)^2}\, -\, \frac{2}{\Gamma(\beta)\Gamma(2\alpha+\beta)}\, <\, 0$$
as soon as $\beta > h(\alpha),$ as is clear from the fact that $\partial_y F(x,y) > 0$ in the proof of Lemma~\ref{lem:h_Function}. By \eqref{LCV} and since $E_{\alpha,\beta}(0) > 0,$ we get
$| E_{\alpha,\beta}(z)| \, <\, E_{\alpha,\beta}(\Re z)$ for $z$ small enough.
\end{proof}

It is plain from the proof of Theorem \ref{thm:Ineq_LE} (b) that if \eqref{eq:Ineq_GE} does not hold globally, then $E_{\alpha,\beta}$ must have non-real zeros. We can hence deduce the following result, which is also relevant for $\alpha \ge 2.$ 
\medskip

\begin{corollary}
\label{ZZ}
The function $E_{\alpha,\beta}$ has non-real zeros whenever $\beta > h(\alpha).$ 
\end{corollary}

Notice that Theorem 3.1.4. in~\cite{PopovSedletski2011} claims that $E_{\alpha,\beta}$ has non-real zeros for $\alpha > 2$ and $\beta > \mu(\alpha) = \alpha^2/\log(2) - \alpha +0.9,$ which is a slightly better result than ours as $\alpha\to\infty$, in view of Remark \ref{ZZZ} (b) with $11\log(2)/3 -1 = 1.54153.. > 0.9.$ On the other hand, our boundary function $h(\alpha)$ is better as $\alpha\to 2$ since $h(2) =4.37228.. < \mu(2) = 4.67078..$ However, this function $h(\alpha)$ it is not optimal either in the neighbourhood of 2. This is illustrated by the following seemingly unnoticed result, which shows that $E_{2,\beta}$ has non real zeroes if and only if $\beta > 3$.

\medskip

\begin{proposition}\label{E24}
{\em (a)} Suppose $\beta\ge 4.$ Then, one has $E_{2,\beta}(x) > 0$ for every $x\in{\mathbb R}.$ In particular, the function  $E_{2,\beta}$ has only non-real zeros and does not satisfy \eqref{eq:Ineq_GE} globally.

\medskip

{\em (b)} Suppose $\beta \in {(3,4)}$. Then, the function $E_{2,\beta}$ has a finite number of negative zeros and does not satisfy \eqref{eq:Ineq_GE} globally.
\end{proposition}

\begin{proof} We begin with the case $\beta = 4.$ Since 

$$E_{2,4}(z) \, =\, \frac{1}{z}\left(\frac{\sinh\sqrt{z}}{\sqrt{z}} -1\right),$$
for every non-zero $z\in\Cset,$ we have

$$E_{2,4}(-x) \, =\, \frac{1}{x}\left(1\, -\, \frac{\sin\sqrt{x}}{\sqrt{x}}\right)\, > \, 0\qquad \mbox{for every $x > 0,$}$$
which clearly implies that $E_{2,4}$ is positive on the negative, hence on the whole, real line. The same property holds now for $E_{2,\beta}$ for all $\beta > 4$ in view of the integral formula (4.4.5) in \cite{GorenfloKilbasMainardiRogosin2020}, which implies (correcting the wrong prefactor therein),

$$E_{2,\beta}(x)\, =\, \frac{1}{\Gamma(\beta-4)}\int_0^1 E_{2,4} (xt^2)\, t^3\, (1-t)^{\beta-5}\, dt\, > \, 0\qquad \mbox{for every $x\in{\mathbb R.}$}$$
For all $\beta \ge 4,$ this implies that $E_{2,\beta}$ has only non-real zeros, and clearly the inequality \eqref{eq:Ineq_GE} cannot hold globally.

\medskip

We next proceed to the case $\beta \in{(3,4)}$. The asymptotic formula 18.1.(22) in \cite{EMOT} yields

$$E_{\alpha,\beta}(-x)\, =\, \frac{1}{x}\, +\, {\mathrm O} \big(x^{(1-\beta)/2}\big)\qquad \mbox{as $x\to\infty,$}$$
showing that there exists $x_\ast > 0$ such that $E_{\alpha,\beta}(-x) > 0$ for all $x \ge x_\ast$ and, by the principle of isolated zeros, that $E_{\alpha,\beta}$ vanishes a finite number of times on the negative half-line. Moreover, we know from Theorem 4.7 in \cite{GorenfloKilbasMainardiRogosin2020} that there exist non-real $z_n\sim -4\pi^2 n^2$ as $n\to\infty$ such that

$$\vert E_{\alpha,\beta} (z_n)\vert \, = \, 0 \, <\, E_{\alpha,\beta} (\Re z_n)$$
for $n$ large enough, so that \eqref{eq:Ineq_GE} does not hold globally either. 
\end{proof}

\begin{remark} 
\label{b34}
{\em (a)} In the case $\beta \in {(3,4)}$, the above proof also shows that the number of negative zeros of $E_{2,\beta}$, counted with their multiplicities, is even. As for $E_{\alpha, 1}$ with $\alpha\in {(1,2)}$, we might think that this number is always non-zero, a question which remains apparently open.

\medskip
{\em (b)} In case $\beta\ge 4,$ the above result shows that the function $1/E_{2,\beta}$ is not {\rm CM}. Indeed, it is well defined on the whole real line and would be analytic on $\Cset$ if it were {\rm CM}, which is clearly not the case since $E_{\alpha,\beta}$ has complex zeros. 
\end{remark} 

 We now show that for all $\alpha > 2$ and $\beta > 2\alpha -1,$ the inequality \eqref{eq:Ineq_GE} is true everywhere except possibly on some compact set. This contrasts with the case $\alpha \in {(1,2]}$, where the failure of \eqref{eq:Ineq_GE} was observed on large, negative values of~$z$. The reason for this difference is that, when $\alpha > 2,$ the leading terms in the expansion of $E_{\alpha,\beta}$ on the negative half-line have an oscillatory character. This non-trivial result advocates the partial relevance of the boundary function $h(\alpha)$ which determines the validity of \eqref{eq:Ineq_GE} in the neighbourhood of zero, albeit in a non-optimal way.

\medskip

\begin{theorem}\label{prop:a > 2}
  If $\alpha>2$ and $\beta > 2\alpha -1,$ there exists a compact set~$K$ such that~\eqref{eq:Ineq_GE} holds for all $z\not\in K.$ 
\end{theorem}

\begin{proof} For simplicity we will set $z = x +{\rm i} y = \rho\eu^{{\rm i} \theta}$ with $x,y \in {\mathbb R}, \rho \ge 0$ and $\theta \in {({-\pi}, \pi]}$. We will let $\rho\to\infty$ and make a discussion according to the values of $x$ and $\theta.$ We begin with the case $x\ge 0$ and use Formula 18.1.(22) in \cite{EMOT}, which implies

$$E_{\alpha,\beta} (z)\, =\, \frac{z^{\frac{1-\beta}{\alpha}}\, \eu^{z^{1/\alpha}}}{\alpha}\, \left( 1 + {\rm o} (1)\right)\qquad\mbox{as $\rho\to\infty.$}$$ 
We first suppose that $\rho\to\infty$ and $x$ remains bounded. Then there exists a constant $c > 0$ such that

$$\left|\frac{E_{\alpha,\beta}(z)}{E_{\alpha,\beta}(x)}\right|\, \ge\, c\, \vert E_{\alpha,\beta}(z) \vert \,\sim\, \frac{c\,\rho^{\frac{1-\beta}{\alpha}}\, \eu^{\cos(\pi/2\alpha) \rho^{1/\alpha}}}{\alpha}\, \to\,\infty.$$
We next suppose $x,\rho\to\infty$ with $\vert\theta\vert\ge \delta \rho^{-1/2\alpha}$ for some $\delta > 0.$ We then have

$$\left|\frac{E_{\alpha,\beta}(z)}{E_{\alpha,\beta}(x)}\right|\, =\, \eu^{f(\rho,\theta)}\, \left( 1 + {\rm o} (1)\right) $$
with
$$f(\rho,\theta)\, =\,\frac{\beta -1}{\alpha}\,\log(\cos\theta)\, +\, \rho^{1/\alpha} (\cos (\theta/\alpha) - (\cos\theta)^{1/\alpha}).$$
The function $\theta\mapsto \cos (\theta/\alpha) - (\cos\theta)^{1/\alpha}$ has derivative $\frac{\sin\theta}{\alpha} ((\cos\theta)^{1/\alpha-1} -1) > 0$ for all $\theta\in {(0,\pi/2)}$, 
and so for $\vert\theta\vert\ge\pi/3$ one has

$$f(\rho,\theta) \, \ge \, \frac{\beta -1}{\alpha}\,\log(\cos\theta)\, +\, \rho^{1/\alpha} (\cos (\pi/3\alpha) - 2^{-1/\alpha})\, \to\, \infty.$$
If $\theta\in[\delta \rho^{-1/2\alpha},\pi/3],$ we compute
\begin{align*}
\partial_\theta f(\rho,\theta) & =  \frac{1}{\alpha} \left( (1-\beta)\tan\theta \, +\, \rho^{1/\alpha} (\tan\theta\, (\cos\theta)^{1/\alpha} - \sin (\theta/\alpha)) \right)\\
& \ge  \frac{1}{\alpha} \left( (1-\beta)\tan\theta \, +\, \rho^{1/\alpha} (\sin\theta - \sin (\theta/\alpha)) \right)\\
& \ge   \frac{1}{\alpha} \left( (1-\beta)\tan\theta \, +\, \rho^{1/\alpha} (1/2 - 1/2\alpha)\,\theta \right)\, > \, 0
\end{align*}
for $\rho$ large enough. By the even character of $\theta\mapsto f(\rho,\theta),$ this implies that for $\delta \rho^{-1/2\alpha}\le \vert\theta\vert\le \pi/3$ one has

$$f(\rho,\theta)\, \ge\,\frac{\beta -1}{\alpha}\,\log(\cos (\delta \rho^{-1/2\alpha}))\, +\, \rho^{1/\alpha} (\cos (\delta \rho^{-1/2\alpha}/\alpha) - (\cos(\delta \rho^{-1/2\alpha}))^{1/\alpha})\, \to\, \frac{\delta^2(\alpha-1)}{2\alpha^2}\, >\, 0,$$
so that 
$$\left|\frac{E_{\alpha,\beta}(z)}{E_{\alpha,\beta}(x)}\right|\, >\, 1$$
for $\rho$ large enough. 
We finally suppose $x,\rho\to\infty$ with $\theta = {\rm o} (\rho^{-1/2\alpha})$.
We will use the expansion

$$\vert E_{\alpha,\beta} (z)\vert^2 \, \sim\, \sum_{k\ge 0}  \frac{(-1)^k y^{2k}}{(2k)!}\, F_k(x)$$
with
$$F_k(x)\, =\, \sum_{i=0}^{2k} \binom{2k}{i} (-1)^i E_{\alpha,\beta}^{(i)} (x) \,E_{\alpha,\beta}^{(2k-i)} (x),$$
which is a direct consequence of~\eqref{Expa}. For $k =1,$ we have seen above that

$$-F_1(x)\, =\, F_{\alpha,\beta} (x)\, \sim\, \frac{(\alpha-1)\, x^{\frac{3-2(\alpha+\beta)}{\alpha}}\, \eu^{2x^{1/\alpha}}}{\alpha^4}\cdot$$
For $k\ge 2,$ we appeal to the decomposition
\begin{equation}
    \label{decompo}
E_{\alpha,\beta}^{(i)}\, =\, \frac{1}{\alpha^i}\sum_{j\ge 0} a_{j,i} E_{\alpha,i(\alpha -1) +\beta + j},    
\end{equation}
where the double sequence \{$a_{j,i}, \, i,j\ge 0\}$ is defined recursively by  $a_{0,i} = 1, a_{j,i} = 0$ for $j > i$ and 
\begin{equation}
    \label{Rec}
    a_{j,i+1}\, - \, a_{j,i} \, =\, (2-\beta + i(1-\alpha) - j)a_{j-1,i}
\end{equation}
 for all $i\ge 0$ and $j = 1,\ldots, i+1,$ which is easily obtained by identifying the coefficients (see also Propositions 2 and 3 in \cite{GP1} for an alternative derivation formulas). In particular, one has 
$$a_{1,i} \, =\, (1-\beta)\,\binom{i}{1}\, + \, (1-\alpha)\, \binom{i}{2}$$
and
$$a_{2,i} \, =\, (1-\beta)(1-\alpha-\beta)\, \binom{i}{2}\, +\, (1-\alpha)(4-2\alpha - 3\beta) \,\binom{i}{3} \, +\, 3(1-\alpha)^2 \,\binom{i}{4}$$
for all $i \ge 0.$ In general, the recursive formula \eqref{Rec} combined with the well-known explicit summation
$$\sum_{k=0}^{i-1} \binom{k}{j}\, =\, \binom{i}{j+1}$$
for all $i,j \ge 1,$ implies that for all $j\ge 0$ there exist bivariate polynomial coefficients $P_{j,j}(\alpha,\beta), \ldots, P_{j,2j}(\alpha,\beta)$ of degree smaller than $j$ such that
$$a_{j,i} \, =\, \sum_{k=j}^{2j} P_{j,k}(\alpha,\beta)\, \binom{i}{k}.$$
In particular, this shows that for all $j\ge 0$ the coefficient $a_{j,i}$ is a polynomial of degree $2j$ in the variate $i$ which vanishes at $i = 0, \ldots, j-1$.  A well-known property of binomial coefficients is that
$$\sum_{i=0}^{2k} \binom{2k}{i} (-1)^i P(i) \, =\, 0$$
for every polynomial $P$ of degree smaller than $2k-1,$ whence
$$\sum_{i=0}^{2k} \binom{2k}{i} (-1)^i a_{j,i}a_{l, 2k-i} \, =\, 0$$
for all $j,l\ge 0$ with $j+l\le k-1.$  A repeated use on \eqref{decompo} of Formula 18.1.(22) in \cite{EMOT}  shows that
$$E_{\alpha,\beta}^{(i)}\, =\, \frac{1}{\alpha^i}\sum_{j\ge 0} a_{j,i} \, x^{\frac{1-\beta - j + i(1-\alpha)}{\alpha}}\, \eu^{x^{1/\alpha}}\, +\, {\rm O} \big( \eu^{c_\alpha x^{1/\alpha}}\big)$$
for some $c_\alpha < 1.$ Putting everything together with Proposition \ref{LCX} implies then
\begin{eqnarray*}
F_k(x) &= & \frac{1}{\alpha^{2k}} \sum_{i,j,l=0}^{2k} \binom{2k}{i} (-1)^i a_{j,i}a_{l, 2k-i}\, x^{\frac{2(1-\beta) + 2k(1-\alpha) - j -l }{\alpha}}\, \eu^{2x^{1/\alpha}}\, +\, {\rm O} \big( \eu^{2c_\alpha x^{1/\alpha}}\big)\\
&= & \frac{1}{\alpha^{2k}} \sum_{j+l\ge k} \left(\sum_{i=0}^{2k}\binom{2k}{i} (-1)^i a_{j,i}a_{l, 2k-i} \right) \, x^{\frac{2(1-\beta) + 2k(1-\alpha) - j -l }{\alpha}}\, \eu^{2x^{1/\alpha}}\, +\, {\rm O} \big( \eu^{2c_\alpha x^{1/\alpha}}\big)  \\
&= &{\rm O} \Big( x^{\frac{k+2 - 2(k\alpha+\beta)}{\alpha}}\, \eu^{2x^{1/\alpha}}\Big) \;  = \; {\rm O} \Big( x^{2(1-k)+\frac{k -1}{\alpha}}\, F_{\alpha,\beta}(x)\Big)
\end{eqnarray*}
for all $k\ge 2$. Therefore, we obtain
\begin{align*}
\vert E_{\alpha,\beta} (z)\vert^2\, -\, \left( E_{\alpha,\beta} (x)\right)^2 & =  y^2 F_{\alpha,\beta} (x)\, \bigg( 1\, +\, \sum_{k\ge 2}  \frac{(-1)^k y^{2(k-1)} \, F_k(x)}{(2k)!\, F_{\alpha,\beta}(x)}\bigg)\\
& =  y^2 F_{\alpha,\beta} (x)\, \bigg( 1\, +\, \sum_{k\ge 2} \; (\tan\theta)^{2(k-1)} \, {\rm O} \big( x^{\frac{k -1}{\alpha}}\big)\bigg)\, \ge\,  \frac{y^2 F_{\alpha,\beta} (x)}{2}\, >\, 0
\end{align*}
for $\rho$ large enough, since $F_{\alpha,\beta}(x)\ge 0$ and $\theta = {\rm o} (\rho^{-1/2\alpha}) = {\rm o} (x^{-1/2\alpha})$. 
Now if~\eqref{eq:Ineq_GE} failed for an unbounded sequence in the right half-plane, then it is easy to see that some subsequence would
fall into one of the asymptotic regimes for the values of~$x,\rho$ we have just covered.
Thus, we have shown that there exists a  compact set~$K_+ \subset \{\Re z\ge 0\}$ such that~\eqref{eq:Ineq_GE} holds for all $z\in \{\Re z\ge 0\}\setminus  K_+$.

\medskip

The argument in the case $x = \Re z\le 0$ is analogous, except that the leading terms of $E_{\alpha,\beta} (z)$ will have an oscillatory character in the vicinity of the negative axis. We first suppose $\rho\to\infty$ while $x$ remains bounded. Then, Formula 18.1.(22) in \cite{EMOT} implies

$$\vert E_{\alpha,\beta} (z)\vert\, \sim\, \frac{\rho^{\frac{1-\beta}{\alpha}}\, \eu^{\cos(\pi/2\alpha)\,\rho^{1/\alpha}}}{\alpha}\; > \; \vert E_{\alpha,\beta} (x)\vert$$
for $\rho$ large enough, since $E_{\alpha,\beta} (x)$ remains bounded, too. We now suppose $\rho\to\infty$ and $x\to -\infty,$ with
\begin{equation}
\label{Oszil}
E_{\alpha,\beta} (x)\, =\, \frac{2}{\alpha}\, \vert x\vert^{\frac{1-\beta}{\alpha}}\, \eu^{\cos(\pi/\alpha)\,\vert x\vert^{1/\alpha}}\cos\big(\vert x\vert^{1/\alpha} \sin(\pi/\alpha) +\pi(1-\beta)/\alpha\big)\left(1\, +\, {\rm o}(1)\right),
\end{equation}
again by Formula 18.1.(22) in \cite{EMOT}. The formula $\cos(x) \cos(y) = \cos^2((x+y)/2) - \sin^2((x-y)/2)$ combined with some algebra leads then to

$$F_{\alpha,\beta} (x)\, \sim\,  \frac{4\, \sin^2(\pi/\alpha)}{\alpha^2}\, \vert x\vert^{\frac{2(2- \alpha-\beta)}{\alpha}}\, \eu^{2\cos(\pi/\alpha)\vert x\vert^{1/\alpha}}\, >\, 0.$$
Similarly as above, we have

$$F_k(x)\, =\, {\rm O} \left( E_{\alpha,\beta} (x) E_{\alpha, 2k(\alpha-1) +\beta +k}(x) \right) \, =\, {\rm O} \left( x^{2(1-k)+\frac{k -2}{\alpha}}\, F_{\alpha,\beta}(x)\right)$$
for all $k\ge 2,$ which by the positivity of $F_{\alpha,\beta} (x)$ implies

$$\vert E_{\alpha,\beta} (z)\vert^2\, -\, \left( E_{\alpha,\beta} (x)\right)^2 \, \ge\, \frac{y^2 F_{\alpha,\beta} (x)}{2}\, >\, 0$$
for $\rho$ large enough, as soon as $\varphi = \pi -\theta = {\rm o} (\rho^{-1/2\alpha}).$ If we next suppose $\vert\varphi\vert \ge \delta > 0$ as $\rho\to\infty,$ we find

$$\left|\frac{E_{\alpha,\beta}(x)}{E_{\alpha,\beta}(z)}\right|\, =\, 2\, \cos(\varphi)^{\frac{1-\beta}{\alpha}}\, \eu^{-f_\alpha(\varphi)\rho^{1/\alpha}}\,\left|\, \cos (\rho^{1/\alpha} (\cos\varphi)^{1/\alpha} \sin(\pi/\alpha) +\pi(1-\beta)/\alpha)\,\right| \left(1\, +\, {\rm o}(1)\right)\, \to\, 0,$$
since $f_\alpha(\varphi) = \cos(\varphi/\alpha) -\cos(\pi/\alpha)\,(\cos\varphi)^{1/\alpha} > \cos(\delta/\alpha) -(\cos\delta)^{1/\alpha} > 0.$ We finally suppose $\varphi = {\rm o}(1)$ with $\vert \varphi\vert \ge \delta \rho^{-1/2\alpha}$ for some $\delta > 0,$ as $\rho\to\infty.$ We may and will suppose $\varphi\ge 0.$ On the one hand, \eqref{Oszil} gives the upper bound

$$\vert E_{\alpha,\beta} (x)\vert\, \le\, \frac{3}{\alpha}\, \rho^{\frac{1-\beta}{\alpha}}\, \eu^{\cos(\pi/\alpha) \,(\cos\varphi)^{1/\alpha}\rho^{1/\alpha}}$$
for $\rho$ large enough. On the other hand, Formula 18.1.(22) in \cite{EMOT} implies

$$E_{\alpha,\beta} (z)\, =\, \frac{1}{\alpha}\, \left( z^{\frac{1-\beta}{\alpha}}\, \eu^{z^{1/\alpha}} \, +\, (z \eu^{-2{\rm i}\pi})^{\frac{1-\beta}{\alpha}}\, \eu^{ (z \eu^{-2{\rm i}\pi})^{1/\alpha}}\right)\,\left(1\, +\, {\rm o}(1)\right)$$
with
\[
  \left|(z \eu^{-2{\rm i}\pi})^{\frac{1-\beta}{\alpha}}\, \eu^{ (z \eu^{-2{\rm i}\pi})^{1/\alpha}}\right|\, = \, \rho^{\frac{1-\beta}{\alpha}}\, \eu^{\cos((\pi+\varphi)/\alpha)\,\rho^{1/\alpha}}\, \ll\, \rho^{\frac{1-\beta}{\alpha}}\, \eu^{\cos((\pi-\varphi)/\alpha)\,\rho^{1/\alpha}}\, =\,\left|z^{\frac{1-\beta}{\alpha}}\, \eu^{z^{1/\alpha}} \right|,
\]
because
\begin{align*}
  (\cos((\pi+\varphi)/\alpha) - \cos((\pi-\varphi)/\alpha))\,\rho^{1/\alpha} &= 2\sin(2\pi/\alpha) \sin(2\varphi/\alpha) \rho^{1/\alpha}\\
  &\ge \delta \sin(2\pi/\alpha) \rho^{1/2\alpha}\,\to\,\infty.
\end{align*}
This finally implies

$$\vert E_{\alpha,\beta} (z)\vert\, \sim\, \frac{\rho^{\frac{1-\beta}{\alpha}}\, \eu^{\cos((\pi-\varphi)/\alpha)\,\rho^{1/\alpha}}}{\alpha}\, \gg\, \vert E_{\alpha,\beta} (x)\vert,$$
since
\begin{align*}
\cos((\pi-\varphi)/\alpha) - \cos(\pi/\alpha) (\cos\varphi)^{1/\alpha} & =  \cos(\pi/\alpha) \left(\cos(\varphi/\alpha) -  (\cos\varphi)^{1/\alpha}\right) + \sin(\pi/\alpha) \sin(\varphi/\alpha)\\
& \ge   \sin(\pi/\alpha) \sin(\varphi/\alpha)   
\end{align*}
and $\sin(\pi/\alpha) \sin(\varphi/\alpha) \rho^{1/\alpha}\,\to\,\infty$. All in all, we have shown that there exists a  compact set~$K_- \subset \{\Re z\le 0\}$ such that~\eqref{eq:Ineq_GE} holds for all $z\in \{\Re z\le 0\}\setminus  K_-$. This completes the proof.
\end{proof}

\begin{remark}
\label{Faul}
It would be interesting to derive an exact formula for the bivariate polynomials $P_{j,k}(\alpha,\beta)$ which were mentioned in the above proof. It is easily shown that 
$$P_{j,j}(\alpha,\beta) \, =\, \prod_{k=0}^{j-1} (1-\beta - k\alpha)\qquad\mbox{and}\qquad P_{j,2j} (\alpha,\beta) \, =\, \frac{(2j)!}{2^j j!}\, (1-\alpha)^j.$$
In general, one has the factorization $P_{j,k}(\alpha,\beta) = (1-\alpha)^{k-j} Q_{j,k}(\alpha,\beta)$ for all $k=j,\ldots, 2j,$ where $Q_{j,k}(\alpha,\beta)$ is a polynomial of degree $2j-k.$ But the combinatorics lying behind these $Q_{j,k}(\alpha,\beta)$ do not seem simple at first sight.
\end{remark}

We finally show that the boundary function $h$ is also particularly relevant for characterizing log-convexity or log-concavity properties of ML functions on the positive half-line. The following result, which builds on Remark 7 in \cite{FerSim}, can be viewed as an extension of Proposition 2.3 in \cite{GerSim}.

\begin{figure}[htb]
\centering
\includegraphics[width=0.6\textwidth]{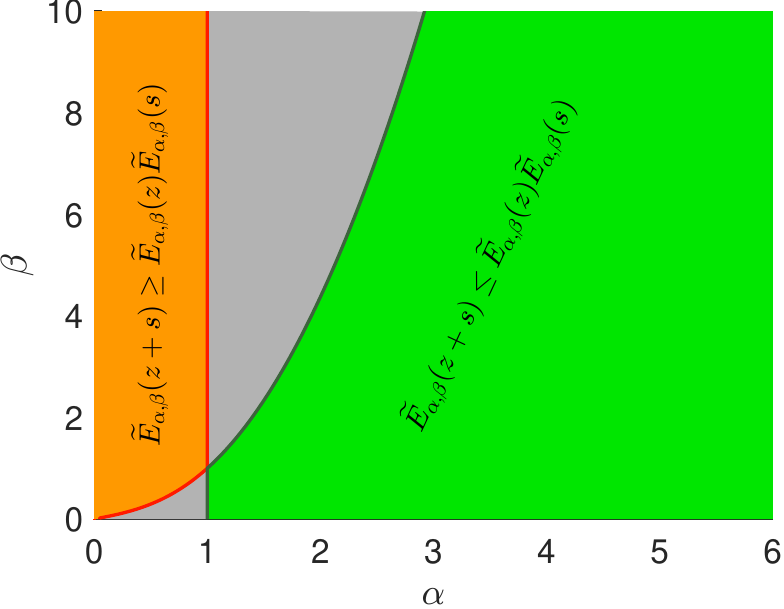}
\caption{
Regions in the parameter space $(\alpha,\beta)$ where $\tilde {E}_{\alpha,\beta}(x)=\Gamma(\beta) E_{\alpha,\beta}(x)$ is super-additive (orange) or sub-additive (green) on $\Rset^{+}.$ The curve represents the separating function $\beta=h(\alpha)$. For $(\alpha,\beta)$ in the gray region, the function $\Gamma(\beta) E_{\alpha,\beta}(x)$ is neither super-additive nor sub-additive on $\Rset^{+}.$}\label{fig:SuperSub}
\end{figure}

\medskip

\begin{theorem} 
\label{LCLX}
Suppose $\alpha > 0$ and $\beta \ge 0.$ Then, one has
$$E_{\alpha, \beta} \;\mbox{is log-convex on ${\mathbb R}^+$}\quad\Longleftrightarrow \quad \alpha\le 1\;\,\mbox{and}\;\,\beta \, \ge\, h(\alpha)$$
and
$$E_{\alpha, \beta} \;\mbox{is log-concave on ${\mathbb R}^+$}\quad\Longleftrightarrow \quad \alpha\ge 1\;\,\mbox{and}\;\,\beta \, \le\, h(\alpha).$$    
\end{theorem}

\begin{proof}
We begin with the characterization of log-convexity and first suppose $\alpha \le 1$. The only if part comes from the definition of the function $h$ since the log-convexity of $E_{\alpha, \beta}$ on ${\mathbb R}^+$ implies that at zero, which amounts to

$$E_{\alpha, \beta}''(0)E_{\alpha, \beta}(0)\, \ge\, E_{\alpha, \beta}'(0)^2\quad\Longleftrightarrow\quad \frac{2}{\Gamma(\beta)\Gamma(\beta+2\alpha)}\, \ge\, \frac{1}{\Gamma(\beta +\alpha)^2}\quad\Longleftrightarrow\quad \beta \ge h(\alpha).$$
We next show the if part, which is clear if $\beta\ge \alpha$ because $E_{\alpha,\beta}$ is then a moment generating function -- see Remark 7 in \cite{FerSim}. The case $1>\alpha > \beta \ge h(\alpha)$ is however less immediate and we will use the same method as in Proposition 2.3 in \cite{GerSim} and Proposition 3 in \cite{FerSim}. Since

$$\frac{E'_{\alpha, \beta}(x)}{E_{\alpha, \beta}(x)}\, =\, \frac{{\displaystyle \sum_{n\ge 0} \frac{(n+1)x^n}{\Gamma(\beta +\alpha +\alpha n)}}}{{\displaystyle \sum_{n\ge 0} \frac{x^n}{\Gamma(\beta +\alpha n)}}}$$
is a ratio of entire series with positive coefficients, it follows by Biernacki and Krzy\.z's lemma -- see \cite{BK} -- that $E_{\alpha, \beta}$ is log-convex on ${\mathbb R}^+$ as soon as the sequence

$$u_n \, =\, \frac{(n+1) \Gamma(\beta +\alpha n)}{\Gamma(\beta +\alpha + \alpha n)}$$
is non-decreasing. The condition $\beta \ge h(\alpha)$ means that $u_1\ge u_0$ and we will show that it also implies $u_{n+1} \ge u_n$ for every $n\ge 1.$ Putting together Frullani's integral and Malmst\'en's formula -- see respectively the formulas 1.7.2.(18) and 1.9.(1) in \cite{EMOT} -- yields, after some simplifications, the Laplace representation

$$\log u_{n+1} - \log u_n\, =\, \int_0^\infty \eu^{-nt} \,\varphi(t)\, dt$$
for every $n\ge 0,$ with

$$\varphi (t) \, =\, \frac{1- \eu^{-t}}{t(1- \eu^{-t/\alpha})}\left(-\eu^{-(\beta/\alpha)t} + \eu^{-t} +\eu^{-(1+\beta/\alpha)t} - \eu^{-(1+1/\alpha)t} \right).$$
The key observation is that the function $\varphi$ strictly changes its sign only once on $(0,\infty),$ starting positive. Indeed, the number of strict sign changes of $\varphi$ on $(0,\infty)$ is that of the function ${\tilde\varphi} (u) = u^{\beta/\alpha} - u - u^{1+\beta/\alpha} + u^{1+1/\alpha}$ on $(0,1)$ and since we have $\beta/\alpha < 1 < 1 +\beta /\alpha < 1+1/\alpha$ for all $\beta\in[h(\alpha), \alpha),$ the generalized Descartes' rule of signs -- see e.g. Theorem 2.2. in \cite{Hauk} -- implies that ${\tilde\varphi}$ vanishes at most twice on $(0,1).$ But we also have ${\tilde\varphi}' (0+) = \infty$ and ${\tilde\varphi}' (1-) = 1/\alpha - 1 >0,$ which implies that ${\tilde\varphi}$ can only have one strict change of sign on $(0,1)$ and so does $\varphi$ on $(0,\infty),$ starting positive because ${\tilde\varphi}' (1-) = 1/\alpha - 1 >0$. Recalling $u_1\ge u_0,$ this shows alltogether that there exists $t_0 > 0$ such that $\varphi(t) > 0$ for $t\in (0,t_0)$ and $\varphi(t) < 0$ for $t\in (t_0, \infty),$ with

$$\int_0^{t_0}\varphi(t)\, dt \, \ge \, -\int_{t_0}^\infty \varphi(t)\, dt\, >\, 0.$$
If we next define positive random variables $X, Y$ with respective densities

$$\frac{\varphi(t)\, {\bf 1}_{(0,t_0)}(t)}{{\displaystyle  \int_0^{t_0}\varphi(t)\, dt}}\qquad\qquad\mbox{and}\qquad\qquad \frac{\varphi(t)\, {\bf 1}_{(t_0,\infty)}(t)}{{\displaystyle \int_{t_0}^\infty\varphi(t)\, dt}},$$
we obviously have ${\mathbb E} [\eu^{-nX}]\ge {\mathbb E} [\eu^{-nY}]$ for all $n\ge 0$ since ${\mathbb P}[X\le Y] = 1.$ This yields 

$$\int_0^{t_0}\eu^{-nt} \, \varphi(t)\, dt\, =\, \left( \int_0^{t_0}\varphi(t)\, dt\right) {\mathbb E} [\eu^{-nX}]  \, \ge \, -\left( \int_{t_0}^\infty \varphi(t)\, dt\right) {\mathbb E} [\eu^{-nY}]\, =\, -\int_{t_0}^\infty \eu^{-nt}\, \varphi(t)\, dt$$
and concludes the argument for the if part. Finally, the estimate $E_{\alpha,\beta} (x) \sim \alpha^{-1} x^{(1-\beta)/\alpha}\eu^{x^{1/\alpha}}$ as $x\to\infty$ implies $\log E_{\alpha, \beta}(x)\sim x^{1/\alpha}$ and shows that $E_{\alpha,\beta}$ cannot be log-convex on ${\mathbb R}^+$ whenever $\alpha > 1.$ 

\medskip

We next proceed to the characterization of log-concavity. The above estimate at infinity shows that we need only to consider the case $\alpha \ge 1.$ Again, the only if part comes from the fact that the log-concavity {of} $E_{\alpha, \beta}$ at zero amounts to

$$E_{\alpha, \beta}''(0)E_{\alpha, \beta}(0)\, \le\, E_{\alpha, \beta}'(0)^2\quad\Longleftrightarrow\quad \frac{2}{\Gamma(\beta)\Gamma(\beta+2\alpha)}\, \le\, \frac{1}{\Gamma(\beta +\alpha)^2}\quad\Longleftrightarrow\quad \beta \le h(\alpha).$$
In order to establish the if part, we need to show that the above sequence $\{u_n\}$ is non-increasing whenever $\beta \le h(\alpha).$ If $\beta \le \alpha$ we have $\beta/\alpha \le 1$ and $\beta/\alpha +(1+1/\alpha) \le 1 + (1 + \beta/\alpha)$,  so that with the notation of Definition A.2 p.12 in \cite{Marsh}, we have the weak majorization 
$$(1,1+\beta/\alpha) \,\prec^W\,(\beta/\alpha, 1+1/\alpha).$$ 
A consequence of Tomi\'c's lemma -- see Proposition 4.B.2. p.157 in \cite{Marsh} -- applied to the convex decreasing function $x\mapsto \eu^{-tx}$ is that the above function $\varphi$ is non-positive on $(0,\infty)$, and $\{u_n\}$ is hence non-increasing. If $\alpha > 1$ and $\beta\in (\alpha, h(\alpha)],$ then we have $1 < \beta/\alpha < 1 +1/\alpha < 1+\beta/\alpha$ and the same reasoning as above shows that $\varphi$ strictly changes its sign only once on $(0,\infty),$ starting negative, and then that   

$$\int_0^\infty \eu^{-nt} \,\varphi(t)\, dt\, \le \, 0$$
for all $n\ge 0,$ as required.
\end{proof}

\begin{remark}
\label{Subsuper}
{\em (a)} As seen during the proof, the log-convexity of $E_{\alpha,\beta}$ holds on the whole ${\mathbb R}$ for $\alpha\le 1$ and $\beta \ge \alpha.$ Observe that log-concavity and log-convexity properties for $E_{\alpha,\beta}$ do not make sense on ${\mathbb R}$ for the other values of $(\alpha, \beta)$ because $E_{\alpha,\beta}$ takes then negative values on ${\mathbb R}^-$.

\medskip
{\em (b)} In case $\beta > 0,$ if we introduce the function ${\tilde E}_{\alpha, \beta} (x) = \Gamma(\beta) E_{\alpha, \beta}(x),$ then we have ${\tilde E}_{\alpha, \beta} (0) = 1$ and the log-convexity of $E_{\alpha, \beta}$ on 
${\mathbb R}^+$ implies the super-additive property 

$${\tilde E}_{\alpha,\beta}(z+s) \ge {\tilde E}_{\alpha,\beta}(z){\tilde E}_{\alpha,\beta}(s)\qquad \mbox{for all $z,s\in {\mathbb R}^+.$}$$ Moreover, a second order expansion at zero shows easily that the log-convexity at zero is a consequence of super-additivity. Thus, we have shown that for all $\beta > 0,$ one has

$$E_{\alpha, \beta} \;\mbox{is log-convex on ${\mathbb R}^+$}\quad\Longleftrightarrow \quad {\tilde E}_{\alpha, \beta} \;\mbox{is super-additive on ${\mathbb R}^+$}\quad\Longleftrightarrow \quad \alpha\le 1\,\;\mbox{and}\,\;\beta \, \ge\, h(\alpha).$$
Similarly, we have

$$E_{\alpha, \beta} \;\mbox{is log-concave on ${\mathbb R}^+$}\quad\Longleftrightarrow \quad {\tilde E}_{\alpha, \beta} \;\mbox{is sub-additive on ${\mathbb R}^+$}\quad\Longleftrightarrow \quad \alpha\ge 1\,\;\mbox{and}\,\;\beta \, \le\, h(\alpha).$$
See Figure \ref{fig:SuperSub} above for an illustration of these two characterizations.
\end{remark}

\section{Discussion and open problems}

\label{S:Conj}

In this last section we summarize our results and formulate some open questions. In Corollary \ref{IFF1} we have characterized the inequality \eqref{eq:Ineq_LE} on the whole $\Cset$. Combined with the main result of \cite{Sc96}, this shows that

$$\mbox{The inequality \eqref{eq:Ineq_LE} holds globally}\quad \Longleftrightarrow\quad 0<\alpha\le 1\; \mbox{and}\; \beta \ge \alpha\quad\Longleftrightarrow\quad E_{\alpha,\beta}(-x)\;\mbox{is CM on $(0,\infty).$}$$
We have observed in Proposition \ref{LCX} that the global inequality \eqref{eq:Ineq_GE} is relevant only for $\alpha\ge 1,$ and Corollary \ref{CHF} implies the following characterization in the case $\alpha = 1:$

$$\mbox{The inequality \eqref{eq:Ineq_GE} holds globally}\quad \Longleftrightarrow\quad \beta \le 1\quad\Longleftrightarrow\quad 1/E_{1,\beta}(x)\;\mbox{is CM on $(0,\infty).$}$$
In the case $\alpha\in {(1,2)}$, it is rather natural in view of the previous equivalences, Proposition~\ref{CM} and Proposition \ref{prop:Ineq_large_LE},  to state the following.

\medskip

\begin{conjecture} 
\label{conj:rhombus}
For every $\alpha\in{(1,2)}$, one has
$$\mbox{The inequality \eqref{eq:Ineq_GE} holds globally}\quad \Longleftrightarrow\quad \beta \in[\alpha-1,\alpha]\quad\Longleftrightarrow\quad 1/E_{\alpha,\beta}(x)\;\mbox{is {\em CM} on $(0,\infty).$}$$
\end{conjecture}

Observe that the two direct implications of this conjecture are true from Proposition~\ref{prop:Ineq_large_LE} resp.\ Proposition~\ref{CM}. Moreover, we have also shown in Corollary \ref{thm:Ineq_GE12} that if $\beta \in[\alpha-1,\alpha]$, then \eqref{eq:Ineq_GE} holds globally except possibly on some compact set of the open left half-plane. However, let us mention that there exists an entire function~$E$ with positive coefficients such that $1/E$ is CM on $(0,\infty)$ and \eqref{eq:Ineq_GE} does not hold globally. For instance\footnote{We are grateful to Mateusz Kwa\'snicki for pointing out this example to us.}, the function

$$\eta_n(z)\, =\, \log \bigg((1+z)(1+z/2)(1+z/4)(1+z/(3+{\rm i}))(1+z/(3-{\rm i})) \prod_{k \ge n} (1+z/(k+{\rm i}))(1+z/(k-{\rm i}))\bigg)$$
has Frullani representation

$$\eta_n(z) \,= \,\int_0^\infty(1- \eu^{-zt}) \,\frac{f_n(t)}{t}\, dt\qquad\mbox{with}\qquad f_n(t) = \eu^{-t} + \eu^{-2t} + \eu^{-4t} + 2 \cos t \bigg(\eu^{-3t} + \sum_{k \ge n} \eu^{-kt}\bigg).$$
It is not difficult to prove that $f_n$ is positive on $(0,\infty)$ for $n$ large enough, which means that~$\eta_n$ is Bernstein, and the entire function 

$$E_n(z)\, =\, \eu^{\eta_n(z)},$$ 
which clearly has positive coefficients, is such that $1/E_n$ is CM on $(0,\infty)$ by Theorem 3.6 in \cite{SSV}. However, one has $E_n({\rm i}-3) = 0 < E_n(-3)$ so that \eqref{eq:Ineq_GE} holds on $\{\Re z \ge 0\}$ but not on the whole complex plane. Notice that in this example, the function $E_n$ has an odd number of negative zeros and an infinite number of complex zeros which are asymptotically close to the real line, whereas the non-real zeros of $E_{\alpha,\beta}$ are asymptotically close to the two half-lines given by $\vert\arg z\vert = \pi\alpha/2$ -- see Theorem~4.7
in~\cite{GorenfloKilbasMainardiRogosin2020}. It seems to us that more structural information on the positive random variable underlying the CM function $1/E_{\alpha,\beta}$ is needed in order to tackle the remaining part of Conjecture~\ref{conj:rhombus}.

\begin{figure}[htb]
\centering
\includegraphics[width=0.6\textwidth]{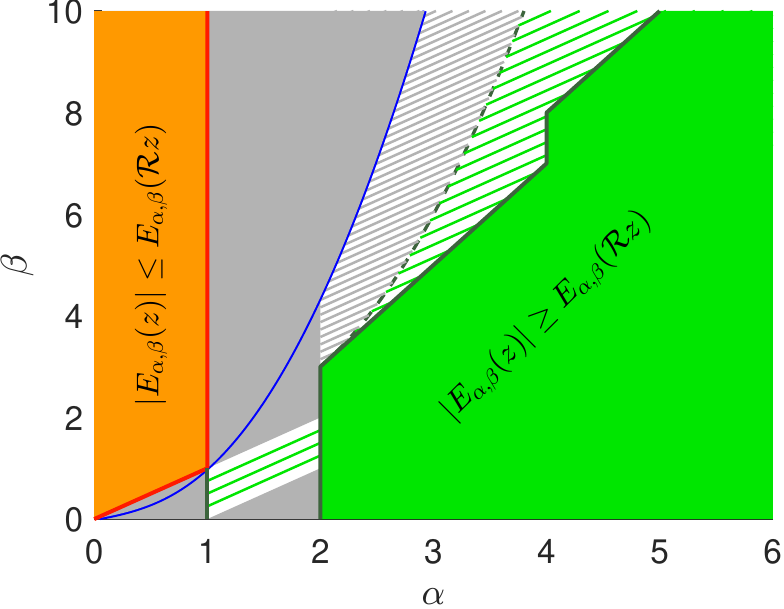}
\caption{
Regions (orange and green) in the parameter space $(\alpha,\beta)$ where inequality \eqref{eq:Ineq_LE} or \eqref{eq:Ineq_GE}, respectively, holds globally, and regions (gray) where neither holds globally. The area with diagonal gray narrow lines corresponds to the region where we conjecture that no inequality holds globally, whereas areas with green spaced diagonal lines correspond to regions where we conjecture that inequality \eqref{eq:Ineq_GE} holds globally. The blue curve represents the function $\beta=h(\alpha)$ from Lemma \ref{lem:h_Function}, while the dashed green curve represents a possible example of the unknown increasing convex function from Conjecture \ref{conj:CBF}.}\label{fig:ParameterSpace_Global}
\end{figure}

In the case $\alpha = 2$, we have obtained the following characterization, as a direct consequence of Theorem \ref{thm:Ineq_LE} (b) and Proposition \ref{E24}:

$$\mbox{The inequality \eqref{eq:Ineq_GE} holds globally}\quad\Longleftrightarrow\quad E_{2,\beta}\;\mbox{has only negative zeros}\quad \Longleftrightarrow\quad \beta \le 3.$$
It is somehow natural to believe that the first equivalence also holds true for $\alpha > 2,$ the reverse implication being already settled in Theorem \ref{thm:Ineq_LE} (b). Proposition \ref{prop:a > 2} gives some partial support for this equivalence, since it shows that for $\alpha > 2$ the inequality \eqref{eq:Ineq_GE} can only fail locally, as is the case for the existence of non-real zeros whose number is necessarily finite. We also believe that the existence of non-real zeros of  $E_{\alpha,\beta}$ for $\alpha > 2$ is characterized by some condition of the type $\beta > f(\alpha),$ as for $\alpha = 2$ with $f(2) = 3.$ This is equivalent to the implication

$$E_{\alpha,\beta}\;\mbox{has non-real zeros}\quad\Longrightarrow\quad E_{\alpha,\beta+\mu}\;\mbox{has non-real zeros}$$
for all $\beta,\mu > 0$ and $\alpha > 2,$ an interesting property of ML functions whose validity does not seem to have been investigated as yet. As for the boundary function $h(\alpha)$ of Lemma \ref{lem:h_Function}, which determines the failure of \eqref{eq:Ineq_GE} in the neighborhood of zero, we may also believe that the boundary function $f(\alpha)$ is increasing and convex. We finally notice that when $E_{\alpha,\beta}$ has only negative zeros, then by the Frullani integral we have the decomposition
\[
  \log\Gamma(\beta)\, +\, \log E_{\alpha,\beta}(z) \, =\, \sum_{n\ge 1} \int_0^\infty(1- \eu^{-z t}) \,\frac{\eu^{z_n t}}{t}\, dt,
\]
where $0 > z_1\ge z_2\ge ...$ are the negative zeros of $E_{\alpha,\beta}.$ This shows that $\log E_{\alpha,\beta}$ is up to translation a complete Bernstein function -- see Chapter~6 in~\cite{SSV} for more details on this notion -- and it is easy to see that this is actually a characterization. We summarize the previous discussion with the following.

\medskip

\begin{conjecture} 
\label{conj:CBF}
There exists an increasing convex function $f : [2,\infty) \to [3,\infty)$ with $f(2) = 3$ such that for every $\alpha \ge 2,$ one has
$$\mbox{The inequality \eqref{eq:Ineq_GE} holds globally}\quad \Longleftrightarrow\quad \log E_{\alpha,\beta}\;\mbox{is a complete Bernstein function}\quad \Longleftrightarrow\quad \beta \le f(\alpha).$$
\end{conjecture}

Observe that if Conjecture \ref{conj:rhombus} is true, then, as seen in Remark~\ref{CBFFF}, we will also have
$$\mbox{The inequality \eqref{eq:Ineq_GE} holds globally}\quad \Longleftrightarrow\quad \log E_{\alpha,\beta}\;\mbox{is a Bernstein function}$$
for all $\alpha\in [1,2).$ However, except in the degenerate case $\alpha =\beta =1,$ Theorem 6.2.(v) in \cite{SSV} implies that the function $\log E_{\alpha,\beta}$ can never be a complete Bernstein function for $\alpha\in [1,2)$, since by the existence of non-real zeros for $E_{\alpha, \beta},$ it does not have an analytic continuation on $\Cset\setminus (-\infty,0].$ All in all, we believe that the characterization of the global inequality \eqref{eq:Ineq_GE} for $\alpha > 1$ is different according as $\alpha < 2$ or $\alpha \ge 2.$ This will be the matter of further research. 

\medskip

The above figure \ref{fig:ParameterSpace_Global} illustrates our current findings and open questions, and concludes the paper.

\section*{Acknowledgments}
We are grateful to an anonymous referee, whose precise and constructive comments helped us improve the quality of the paper.
The work of RG and MP is partially supported by the INdAM
under the GNCS Project CUP E53C24001950001. RG received financial support also by the MUR under the PRIN-PNRR2022 project n.\ P2022M7JZW, grant No.\ CUP H53D23008930001.
SG  received financial support from the Austrian Science Fund (FWF) [10.55776/PAT1474824].
MP received financial support also by the MUR under the PRIN-PNRR2022 project n.\
P202254HT8, grant No.\ CUP B53D23027760001


\end{document}